\documentclass[12pt]{amsart}
\usepackage{amssymb,amsmath}
\usepackage{epsfig,multirow}
\usepackage{graphicx}
\usepackage{colordvi}
\copyrightinfo{}{}

\newtheorem{ithm}{Theorem}

\numberwithin{equation}{section}
\newtheorem{theorem}[equation]{Theorem}

\newtheorem{corollary}[equation]{Corollary}
\newtheorem{conjecture}[equation]{Conjecture}

\theoremstyle{definition}
\newtheorem{definition}[equation]{Definition}

\newtheorem{example}[equation]{Example}

\newcounter{FNC}[page]
\def\fauxfootnote#1{{\addtocounter{FNC}{2}$^\fnsymbol{FNC}$%
     \let\thefootnote\relax\footnotetext{$^\fnsymbol{FNC}$#1}}}

\newcommand{\maxp}{\includegraphics{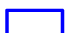}}
\newcommand{\SYT}{\mbox{\rm SYT}}

\newcommand{\C}{{\mathbb C}}
\newcommand{\R}{{\mathbb R}}
\newcommand{\PP}{{\mathbb P}}

\newcommand{\calR}{{\mathcal{R}}}
\newcommand{\calS}{{\mathcal{S}}}
\newcommand{\calB}{{\mathcal{B}}}

\newcommand{\G}{{\mathbb G}}
\newcommand{\Gr}{{\G(n,d)}}

\newcommand{\Fl}{{\mathbb F}\ell_{\bn,d}}
\newcommand{\Flv}{{\mathbb F}\ell}

\newcommand{\gln}{{\mathfrak{gl}_{n+1}}}
\newcommand{\sln}{{\mathfrak{sl}_{n+1}}}
\newcommand{\frakh}{{\mathfrak{h}}}
\newcommand{\frakn}{{\mathfrak{n}}}

\newcommand{\Wr}{{\rm Wr}}
\newcommand{\sing}{{\rm sing}}
\newcommand{\dlog}{{\rm ln}'}
\newcommand{\ord}{{\rm ord}}
\newcommand{\diag}{{\rm diag}}
\newcommand{\slide}{{\rm slide}}

\newcommand{\bA}{{\bf A}}
\newcommand{\ba}{{\bf a}}
\newcommand{\sba}{{\mbox{\scriptsize${\bf a}$}}}
\newcommand{\bb}{{\bf b}}

\newcommand{\bk}{{\mbox{\boldmath$\kappa$}}}
\newcommand{\sbk}{{\mbox{\scriptsize\boldmath$\kappa$}}}
\newcommand{\bmu}{{\mbox{\boldmath$\mu$}}}
\newcommand{\sbmu}{{\mbox{\scriptsize\boldmath$\mu$}}}
\newcommand{\bnu}{{\mbox{\boldmath$\nu$}}}

\newcommand{\blambda}{{\mbox{\boldmath$\lambda$}}}
\newcommand{\sbl}{{\mbox{\scriptsize\boldmath$\lambda$}}}
\newcommand{\bM}{{\bf M}}
\newcommand{\bn}{{\bf n}}
\newcommand{\bp}{{\bf p}}
\newcommand{\bs}{{\bf s}}

\newcommand{\bt}{{\bf t}}
\newcommand{\bu}{{\bf u}}
\newcommand{\bx}{{\bf x}}

\newcommand{\bz}{{\bf z}}

\newcommand{\bi}{{\mbox{\boldmath$\iota$}}}
\newcommand{\sbi}{{\mbox{\scriptsize\boldmath$\iota$}}}

\newcommand{\DeCo}[1]{\Blue{#1}}

\begin{document}
\title{Frontiers of reality in Schubert calculus}
\author{Frank Sottile}
\address{Department of Mathematics\\
         Texas A\&M University\\
         College Station\\
         TX \ 77843\\
         USA}
\email{sottile@math.tamu.edu}
\urladdr{www.math.tamu.edu/\~{}sottile}
\thanks{Work of Sottile supported by NSF grant DMS-0701050}
\begin{abstract}
 The theorem of Mukhin, Tarasov, and Varchenko (formerly the Sha\-pi\-ro conjecture for
 Grassmannians) asserts that all ({\it a priori} complex) solutions to certain
 geometric problems in the Schubert calculus are actually real. Their proof is quite
 remarkable, using ideas from integrable systems, Fuchsian differential equations, and
 representation theory.  
 There is now a second proof of this result, and it has
 ramifications in other areas of mathematics, from curves to control theory to combinatorics.
 Despite this work, the original Shapiro conjecture is not yet settled.
 While it is false as stated, it has several interesting and not quite understood
 modifications and generalizations that are likely true, and the strongest and most subtle
 version of the Shapiro conjecture for Grassmannians remains open. 
\end{abstract}
\maketitle

\section*{Introduction}

While it is not unusual for a univariate polynomial $f$ with real coefficients to 
have {\sl some}
real roots---under reasonable assumptions we expect $\sqrt{\deg f}$ real
roots~\cite{Ko93}---it is rare for a polynomial to have all of its roots real.
In a sense, the only natural example of a polynomial with all of its roots real is the
characteristic polynomial of a real symmetric matrix, as all eigenvalues of a symmetric matrix
are real. 

Similarly, when a system of real polynomial equations has finitely many
({\it a priori} complex) solutions, we expect some, but likely not all, solutions to be real.
In fact, upper bounds on the number of real solutions~\cite{BBS,Kh80} sometimes ensure that not
all solutions can be real. 
As before, the most natural example 
of a system with only real 
solutions is the system of equations for
the eigenvectors and eigenvalues of a real symmetric matrix.

Here is another system of polynomial equations that has only
real solutions.
The Wronskian of univariate polynomials $f_0,\dotsc,f_n\in\C[t]$ is the
determinant
\[
   \det\ \left(\begin{matrix}
        f_0(t) &f_1(t)&\dotsb &f_n(t) \\
        f'_0(t) &f'_1(t)&\dotsb &f'_n(t) \\
         \vdots & \vdots & \ddots & \vdots\\
     f^{(n)}_0(t) &f^{(n)}_1(t)&\dotsb &f^{(n)}_n(t) 
    \end{matrix}\right)\ .
\]
Up to a scalar multiple, the Wronskian depends only upon the
linear span $P$ of the polynomials $f_0,\dotsc,f_n$. 
This scaling retains only the information of the roots and their
multiplicities. 
Recently, Mukhin, Tarasov, and Varchenko~\cite{MTV_Sh} proved the remarkable (but
seemingly innocuous) result.

\begin{ithm}\label{Th:MTV_1}
 If the Wronskian of a vector space $P$ of polynomials has only real roots, then $P$ has a
 basis of real polynomials. 
\end{ithm}

While not immediately apparent, those $(n{+}1)$-dimensional subspaces $P$ of $\C[t]$ with
a given Wronskian $W$ are the solutions to a system of polynomial equations that depend
on the roots of $W$.
In Section~\ref{S:Shapiro}, we explain how the Shapiro conjecture for Grassmannians
is equivalent to Theorem~\ref{Th:MTV_1}.

The proof of Theorem~\ref{Th:MTV_1} uses the Bethe ansatz for the (periodic) Gaudin model on
certain modules (representations) of the Lie algebra $\sln\C$.
The Bethe ansatz is a method to find pure states, called \DeCo{{\sl Bethe vectors}}, of quantum
integrable systems~\cite{Ga}. 
Here, that means common eigenvectors for a family of commuting operators called the
Gaudin Hamiltonians which generate a commutative Bethe algebra $\calB$. 
As $\calB$ commutes with the action of $\sln\C$, this also 
decomposes a module of $\sln\C$ into irreducible submodules.
It includes a set-theoretic map from the  Bethe eigenvectors to 
spaces of polynomials with a given Wronskian. 
A coincidence of numbers, from the Schubert calculus and from representation
theory, implies that this map is a bijection.
As the Gaudin Hamiltonians are symmetric with respect to the positive definite Shapovalov
form, their eigenvectors and eigenvalues are real. 
Theorem~\ref{Th:MTV_1} follows as eigenvectors with real eigenvalues must come from real
spaces of polynomials. 
We describe this in Sections~\ref{S:polys}, \ref{S:BAGM}, and~\ref{S:Shapovalov}.

There is now a second proof~\cite{MTV_R} of Theorem~\ref{Th:MTV_1}, also passing
through integrable systems and representation theory. 
It provides a deep connection
between the Schubert calculus and the representation theory of $\sln\C$, 
strengthening Theorem~\ref{Th:MTV_1} to include transversality.

The geometry behind the statement of Theorem~\ref{Th:MTV_1} appears in many other guises,
some of which we describe in Section~\ref{S:applications}.
These include linear series on the projective line~\cite{EH83}, 
rational curves with prescribed flexes~\cite{KS}, and the feedback control of a system of
linear differential equations~\cite{By89,EG02b}.
A special case of the Shapiro conjecture concerns 
rational functions with prescribed critical points, and was proved in this form by Eremenko and 
Gabrielov~\cite{EG02a}.
They showed that a rational function whose critical points lie on a
circle in the Riemann sphere maps that circle to another circle.
Using the strengthening of Theorem~\ref{Th:MTV_1} involving transversality,
Purbhoo~\cite{Purbhoo} discovered that the fundamental combinatorial algorithms on Young
Tableaux come from the monodromy of the map that takes spaces of polynomials to their
Wronskians.  

A generalization of Theorem~\ref{Th:MTV_1} by Mukhin, Tarasov, and
Varchenko~\cite{MTV_XXX} implies the following
attractive statement from matrix theory.
Let $b_0,b_1,\dotsc,b_n$ be distinct real numbers, $\alpha_0,\dotsc,\alpha_n$ be
complex numbers, and consider the matrix
\[
    Z\ :=\ \left(\begin{matrix}
      \alpha_0& (b_0-b_1)^{-1}& \dotsb & (b_0-b_n)^{-1}\\
      (b_1-b_0)^{-1}& \alpha_1& \dotsb & (b_1-b_n)^{-1}\\
         \vdots      & \vdots & \ddots & \vdots\\
      (b_n-b_0)^{-1}& (b_n-b_1)^{-1}& \dotsb & \alpha_n
    \end{matrix}\right)\ .
\]

\begin{ithm}\label{Th:matrix}
  If $Z$ has only real eigenvalues, then $\alpha_1,\dotsc,\alpha_n$ are real.
\end{ithm}

Unlike its proof, the statement of Theorem~\ref{Th:matrix} has nothing to do with Schubert
calculus or representations of $\sln\C$ or integrable systems, and it remains a challenge
to prove it directly. 
We discuss this in Section~\ref{S:other}.

The statement and proof of Theorem~\ref{Th:MTV_1} is only part of this story.
Theorem~\ref{Th:MTV_1} settles (for Grassmannians) a conjecture in Schubert calculus
made by Boris Shapiro and Michael Shapiro in 1993/4.
While this Shapiro conjecture is false for most other flag manifolds, there are appealing
corrections and generalizations supported by theoretical evidence and by overwhelming
computational evidence, and the strongest and most subtle form remains open.
We sketch this in Section~\ref{S:extensions}.

%
%
\subsection*{First steps:  the problem of four lines}
We close this Introduction by illustrating the Schubert calculus and the Shapiro conjecture 
with some beautiful geometry.
Consider the set of all lines in three-dimensional space.
This set (a Grassmannian) is four-dimensional, which we may see by counting the degrees of
freedom for a line $\ell$ as follows.
Fix planes $\Pi$ and $\Pi'$ that meet $\ell$ in points $p$ and $p'$ as shown.
\[
 \begin{picture}(180,120)
   \put(0,0){\includegraphics[height=120pt]{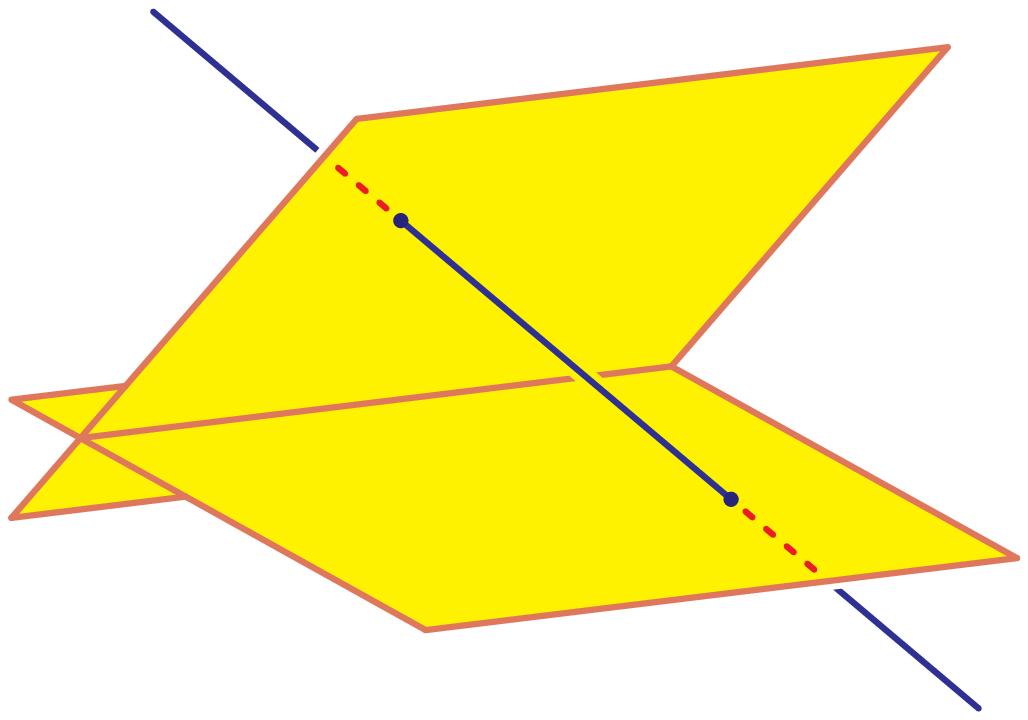}}
   \put(22,105){$\ell$}
   \put(58,78){$p$}   \put(109,29){$p'$}
   \put(155,95){$\Pi$}  \put(165,35){$\Pi'$}
 \end{picture}
\]
Since each point $p, p'$ has two degrees of freedom to move within its 
plane, we see that the line $\ell$ enjoys four degrees of freedom.

Similarly, the set of lines that meet a fixed line is three-dimensional.
More parameter counting tells us that if we fix four lines, then the set of lines that
meet each of our fixed lines will be zero-dimensional.
That is, it consists of finitely many lines.
The Schubert calculus gives algorithms to determine this number of lines.
We instead use elementary geometry to show that this number is 2.

The Shapiro conjecture asserts that if the four fixed lines are chosen in a particular way,
then both solution lines will be real.
This special choice begins by specifying a twisted cubic curve, $\gamma$.
While any twisted cubic will do, we'll take the one with parametrization
\[ 
   \gamma\ \colon\ t\ \longmapsto 
     (6t^2-1, \ \tfrac{7}{2}t^3+\tfrac{3}{2}t,\ \tfrac{3}{2}t-\tfrac{1}{2}t^3)\,.
 \leqno(1) \label{Eq:Twisted}
%
%
\]  
Our fixed lines will be four lines tangent to $\gamma$.

\begin{figure}[htb]
\[
  \begin{picture}(348,220)
   \put(  0,  0){\includegraphics[width=340pt]{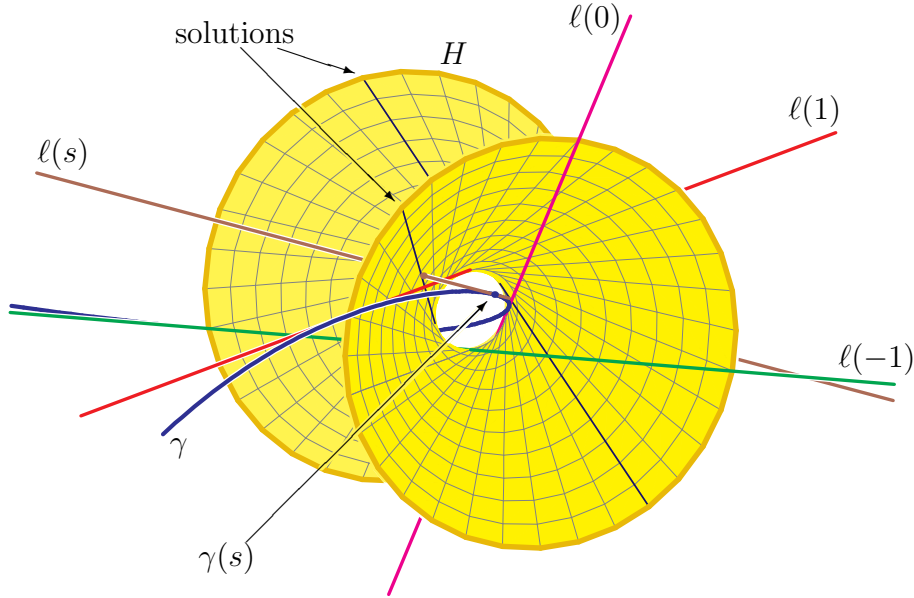}}
   \put( 13,165){$\ell(s)$}
   \put(297,181){$\ell(1)$}
   \put(316, 87){$\ell(-1)$}
   \put(214,216){$\ell(0)$}
   \put( 74,  12){$\gamma(s)$} \put(92,21){\vector(1,1){91}} 
   \put( 63, 55){$\gamma$}
   \put(165,202){$H$}
   \put(65,210){solutions} 
   \put(104,208){\vector(3,-1){29}}
   \put( 90,208){\vector(1,-1){58}}
  \end{picture}
\]
\caption{The problem of four lines.}
\label{F:FRSC}
\end{figure}

We understand the lines that meet our four tangent lines by first considering
lines that meet three tangent lines.
We are free to fix the first three points of tangency to be any of our choosing,
for instance, $\gamma(-1)$, $\gamma(0)$, and $\gamma(1)$.
Then the three lines $\ell(-1)$, $\ell(0)$, and $\ell(1)$  tangent at these points 
have parametrizations
\[
  (-5+s,5-s,-1)\,,\ \ (-1,s,s)\,,\ \ \mbox{\rm and}\ \ 
   (5+s,5+s,1)\ \ \mbox{\rm for $s\in\R$.}
\]
These lines all lie on the hyperboloid $H$ of one sheet defined by
%
\[  \label{Eq:HypEq}
   x^2-y^2+z^2\ =\ 1\,, \leqno(2) 
\]  
which has two rulings by families of lines.
The lines  $\ell(-1)$, $\ell(0)$, and $\ell(1)$ lie in one family, and the other family
consists of the lines meeting $\ell(-1)$, $\ell(0)$, and $\ell(1)$.
This family is drawn on the hyperboloid $H$ in Figure~\ref{F:FRSC}.

The lines that meet $\ell(-1)$, $\ell(0)$, $\ell(1)$, and a fourth line $\ell(s)$ 
will be those in this second family that also meet $\ell(s)$.
In general, there will be two such lines, one for each point of intersection of line
$\ell(s)$ with $H$, as $H$ is defined by the quadratic polynomial~(2).
The remarkable geometric fact 
is that every such tangent line, $\ell(s)$ for $s\not\in\{-1,0,1\}$, will meet the
hyperboloid in two real points. 
We illustrate this when $s=0.31$ in Figure~\ref{F:FRSC}, highlighting the two solution
lines. 

The Shapiro conjecture and its extensions claim that this reality always
happens:
If the conditions for a Schubert problem are chosen in a particular way relative to a
rational normal curve (here, tangent lines to the twisted cubic curve $\gamma$ of~(1)),
then all solutions will be real.  
When the Schubert problem comes from a Grassmannian (like this problem of four lines),
the Shapiro conjecture is true---this is the theorem of Mukhin, Tarasov, and Varchenko.
For most other flag manifolds, it is known to fail, but in very interesting ways.

%
\subsection*{Acknowledgments}
We thank those who have helped us to understand this story and to improve this
exposition.
In particular, we thank Eugene Mukhin, Alexander Varchenko, Milen Yakimov, Aaron Lauve, Zach
Teitler, and Nickolas Hein.

%
%
\section{The Shapiro conjecture for Grassmannians}\label{S:Shapiro}

Let \DeCo{$\C_d[t]$} be the set of complex polynomials of degree at most $d$ in the
indeterminate $t$, a vector space of dimension $d{+}1$.
Fix a positive integer $n<d$ and let \DeCo{$\Gr$} 
be the set of all $({n{+}1})$-dimensional
linear subspaces $P$ of $\C_d[t]$.
This \DeCo{{\sl Grassmannian}} is a complex manifold of dimension 
$(n{+}1)(d{-}n)$~\cite[Ch.~1.5]{GH78}.

The main character in our story is the Wronski map, which associates to a point
$P\in \Gr$ the Wronskian of a basis for $P$.
If $\{f_0(t),\dotsc,f_n(t)\}$ is a basis for $P$, its Wronskian is the determinant of the
derivatives of the basis,
 \begin{equation}\label{Eq:Wronski}
  \DeCo{\Wr}(f_0,\dotsc,f_n)\ :=\ \det\; \left(
  \begin{matrix}
    f_0&f'_0&\dotsb&f^{(n)}_0\\
    f_1&f'_1&\dotsb&f^{(n)}_1\\
    \vdots&\vdots&\ddots&\vdots\\
    f_n&f'_n&\dotsb&f^{(n)}_n\\
  \end{matrix}\right)\ ,
 \end{equation}
which is a nonzero polynomial of degree at most $(n{+}1)(d{-}n)$.
This does not quite define a map $\Gr\to \C_{(n{+}1)(d{-}n)}[t]$, as 
choosing a different basis for $P$ multiplies the Wronskian by
a nonzero constant.
If we consider the Wronskian up to a nonzero constant, we obtain the 
\DeCo{{\sl Wronski map}} 
 \begin{equation}\label{Eq:Wronski_map}
   \DeCo{\Wr}\ \colon\ \Gr\ \longrightarrow\ 
    \PP(\C_{(n{+}1)(d{-}n)}[t])\ \simeq\ \PP^{(n{+}1)(d{-}n)}\,,
 \end{equation}
where $\PP(V)$ denotes the projective space consisting of all $1$-dimensional linear
subspaces of a vector space $V$.

We restate Theorem~\ref{Th:MTV_1}, the simplest version of the Theorem of Mukhin, Tarasov,
and Varchenko~\cite{MTV_Sh}. 
\medskip

\noindent{\bf Theorem~\ref{Th:MTV_1}.}
 {\it
    If the Wronskian of a space $P$ of polynomials has only real roots, then $P$ has a basis
    of real polynomials. 
}\medskip

The problem of four lines in the Introduction is a special case of Theorem~\ref{Th:MTV_1} when
$d=3$ and $n=1$.
To see this, note that if we apply an affine function $a+bx+cy+dz$ to the curve
$\gamma(t)$ of~(1), we obtain a cubic polynomial in $\C_3[t]$, and every
cubic polynomial comes from a unique affine function. 
A line $\ell$ in $\C^3$ (actually in $\PP^3$) is cut out by a two-dimensional space of 
affine functions, which gives a 2-dimensional space $\DeCo{P_\ell}$ of polynomials in
$\C_3[t]$, and hence a point $P_\ell\in\G(1,3)$.

It turns out that the Wronskian point $P_\ell\in\G(1,3)$ is a quartic polynomial with
a root at $s\in\C$ if and only if the corresponding line $\ell$ meets the
line $\ell(s)$ tangent to the curve $\gamma$ at $\gamma(s)$.
Thus a line $\ell$ meets four lines tangent to $\gamma$ at real points if and only if 
the Wronskian of $P_\ell\in \G(1,3)$ vanishes at these four points.
Since these points are real, Theorem~\ref{Th:MTV_1} implies that $P_\ell$ has a basis of real 
polynomials.
Thus $\ell$ is cut out by real affine functions, and hence is real.

%
%
\subsection{Geometric form of the Shapiro conjecture}\label{S:geometry}

Let $P\in\Gr$ be a subspace.
We consider the order of vanishing at a point $s\in\C$ of polynomials in a basis for $P$. 
There will be a minimal order $a_0$ of vanishing for these polynomials.
Suppose that $f_0$ vanishes to this order.
Subtracting an appropriate multiple of $f_0$ from each of the other polynomials,
we may assume that they vanish to order greater than $a_0$ at $s$.
Let $a_1$ be the minimal order of vanishing at $s$ of these remaining polynomials.
Continuing in this fashion, we obtain a basis $f_0,\dotsc,f_n$ of $P$
and a sequence 
\[
    0\leq a_0<a_1<\dotsb<a_n\leq d\,,
\]
where $f_i$ vanishes to order $a_i$ at $s$. 
Call this sequence \DeCo{$\ba_P(s)$} the \DeCo{{\sl ramification}} of $P$ at $s$.
For a sequence $\ba:0\leq a_0<\dotsb<a_n\leq d$, write $\Omega^\circ_\ba(s)$ for the
set of points $P\in \Gr$ with $\ba_P(s)=\ba$, which is a Schubert cell of $\Gr$.
It has codimension 
\[
   \DeCo{|\ba|}\ :=\ a_0\ +\ a_1{-}1\ +\ \dotsb\ +\ a_n{-}n\,,
\]
as may be seen by expressing the basis $f_0,\dotsc,f_n$ of $P$ in terms of the basis
$\{(t-s)^i\mid i=0,\dotsc,d\}$ of $\C_d[t]$.
Since $f^{(i)}_j$ vanishes to order at least $a_j-i$ at $s$ and $f^{(i)}_i$ vanishes to
order exactly $a_i-i$ at $s$, the Wronskian of a subspace
$P\in\Omega^\circ_\ba(s)$ vanishes to order exactly $|\ba|$ at $s$. 

Let \DeCo{$\Gr^\circ$} be the dense open subset of $\Gr$ consisting of those $P$
having a basis $f_0,\dotsc,f_n$ where $f_i$ has degree $d{-}n{+}i$.
When $P\in\Gr^\circ$, we obtain the Pl\"ucker formula for the total ramification of a
general subspace $P$ of $\C_d[t]$,
 \begin{equation}\label{Eq:Pl}
   \dim \Gr\ =\   \sum_{s\in\C} |\ba_P(s)|\,.
 \end{equation}
In general, the total ramification of $P$ is bounded by the dimension of $\Gr$.
(One may also define ramification at infinity for subspaces $P\not\in\Gr^\circ$ to 
obtain the Pl\"ucker formula in its full generality.)
If $\ba_P(s)\colon 0<1<\dotsb<n$, so that $|\ba_P(s)|=0$, then $P$ is 
\DeCo{{\sl unramified}} at $s$.
Theorem~\ref{Th:MTV_1} states that if a subspace $P\in\Gr$ is ramified
only at real points, then $P$ has a basis of real polynomials.

We introduce some more geometry.
Let $W=\prod_s(t-s)^{|\ba_P(s)|}$ be the Wronskian of $P$.
Then
\[
    P\ \in\ \bigcap_{s\colon W(s)=0}\; \Omega_{\ba_P(s)}^\circ (s)\,,
\]
and this intersection consists of all subspaces with the same ramification as $P$.
In particular, $P$ lies in the intersection of the closures of these Schubert cells, 
which we now describe.
For each $s\in\C$, $\C_d[t]$ has a complete flag of subspaces
\[
   \DeCo{F_\bullet(s)}\ :\ 
    \C\cdot (t{-}s)^d\ \subset\ \C_1[t]\cdot (t{-}s)^{d-1}\ \subset\ 
     \dotsb\ \subset\ \C_{d-1}[t]\cdot(t{-}s)\ \subset\ \C_d[t]\,.
\]
More generally, a flag $F_\bullet$ is a sequence of subspaces
\[
   F_\bullet\ :\ F_1\ \subset\ F_2\ \subset\ \dotsb\ \subset\ F_{d}\ \subset \C_d[t]\,,
\]
where $F_i$ has dimension $i$.
For a sequence $\ba$ and a flag $F_\bullet$, the \DeCo{{\sl Schubert variety}} 
 \begin{equation}\label{Eq:Schubert_variety}
   \{ P\in \Gr\mid \dim \left(P\cap F_{d+1-a_j}\right)\geq n{+}1{-}j,\ 
          \mbox{for}\ j=0,1,\dotsc,n\}\,,
 \end{equation}
%
%
is a subvariety of $\Gr$, written $\DeCo{\Omega_{\ba}F_\bullet}$.
It consists of linear subspaces $P$ having special position (encoded by $\ba$) with
respect to the flag $F_\bullet$.
Since $\dim (P\cap F_{d+1-i}(s))$ counts the number of linearly independent polynomials in
$P$ that vanish to order at least $i$ at $s$, we see that 
$\Omega^\circ_\ba(s)\subset\Omega_\ba F_\bullet(s)$.
More precisely, $\Omega_{\ba}F_\bullet(s)$ is the closure of the Schubert cell
$\Omega_{\ba}^\circ(s)$ and it is the disjoint union of cells $\Omega_{\bb}^\circ(s)$ for 
$\bb\geq\ba$, where $\geq$ is componentwise comparison.

Given sequences $\ba^{(1)},\dotsc,\ba^{(m)}$ and flags
$F_\bullet^{(1)},\dotsc,F_\bullet^{(m)}$, the intersection 
 \begin{equation}\label{Eq:Schubert_intersection}
   \Omega_{\ba^{(1)}}F_\bullet^{(1)}\ \bigcap\ 
   \Omega_{\ba^{(2)}}F_\bullet^{(2)}\ \bigcap\ \dotsb\ \bigcap\ 
   \Omega_{\ba^{(m)}}F_\bullet^{(m)}
 \end{equation}
consists of those linear subspaces $P\in G$ having specified position $\ba^{(i)}$ with
respect to the flag $F^{(i)}_\bullet$, for each $i=1,\dotsc,m$. 
Kleiman~\cite{Kl} showed that if the flags $F^{(i)}_\bullet$ are general, then the
intersection~\eqref{Eq:Schubert_intersection} is (generically) transverse.

A \DeCo{{\sl Schubert problem}} is a list $\bA:=(\ba^{(1)},\dotsc,\ba^{(m)})$ of sequences
satisfying
\[
   |\ba^{(1)}|+\dotsb+|\ba^{(m)}|\ =\ (n{+}1)(d{-}n)\ (\,=\ \dim \Gr\,)\,.
\]
Given a Schubert problem, Kleiman's Theorem implies that a general
intersection~\eqref{Eq:Schubert_intersection} will be zero-dimensional and thus consist of
finitely many points.
By transversality, the number \DeCo{$\delta(\bA)$} of these points is
independent of choice of general flags.
The Schubert calculus~\cite{KlLa}, through the Littlewood-Richardson rule~\cite{Fu}, gives
algorithms to determine $\delta(\bA)$.

We mention an important special case.
Let  $\DeCo{\bi}\colon 0<1<\dotsb<n{-}1<n{+}1$ be the unique ramification sequence with
$|\bi|=1$, and write  $\DeCo{\bi_{n,d}}$ for the Schubert problem in which 
$\bi$ occurs $(n{+}1)(d{-}n)$ times.
Schubert~\cite{Sch1886c} gave the formula
 \begin{equation}\label{Eq:Grass_degree}
   \delta(\bi_{n,d})\ = \ [(n{+}1)(d{-}n)]!
          \frac{1!2!\dotsb n!}{(d{-}n)!(d{-}n{+}1)!\dotsb d!}\ .
 \end{equation}

By the Pl\"ucker Formula~\eqref{Eq:Pl}, the total ramification $(\ba_P(s)\colon|\ba_P(s)|>0)$ of
a subspace $P\in\Gr^\circ$ is a Schubert problem.
Let $W$ be the Wronskian of $P$.
We would like the intersection containing $P$,
 \begin{equation}\label{Eq:Shapiro_intersection}
    \bigcap_{s\colon W(s)=0}\; \Omega_{\ba_P(s)} F_\bullet(s)\,,
 \end{equation}
to be transverse and zero-dimensional.
However, Kleiman's Theorem does not apply, as the flags $F_\bullet(s)$ for $s$ a root of
$W$ are not generic. 
For example, in the problem of four lines, if the Wronskian is $t^4-t$, then the
corresponding intersection~\eqref{Eq:Shapiro_intersection} of Schubert varieties is not
transverse. 
(This has been worked out in detail in~\cite[\S 9]{EH83}.)

We can see that this intersection~\eqref{Eq:Shapiro_intersection} is however always
zero-dimensional. 
Note that any positive-dimensional subvariety meets $\Omega_\sbi F_\bullet$, for any flag
$F_\bullet$.
(This is because, for example,  $\Omega_\sbi F_\bullet$ is a hyperplane section of $\Gr$ in
its Pl\"ucker embedding into projective space.)
In particular, if the intersection~\eqref{Eq:Shapiro_intersection} is not zero-dimensional,
then given a point $s\in\PP^1$ with $W(s)\neq 0$, there will be a point $P'$ 
in~\eqref{Eq:Shapiro_intersection}  which also lies in $\Omega_\iota F_\bullet(s)$.
But then the total ramification of $P'$ does not satisfy the Pl\"ucker
formula~\eqref{Eq:Pl}, as its ramification strictly contains
the total ramification of $P$.

A consequence of this argument is that the Wronski map~\eqref{Eq:Wronski_map} is a
flat, finite map.
In particular, it has finite fibers.
The intersection number $\delta(\bi_{n,d})$ in~\eqref{Eq:Grass_degree} is an upper
bound for the cardinality of a fiber.
By Sard's Theorem, this upper bound is obtained for generic Wronskians.
An argument that proves this in somewhat greater generality was given by Eisenbud and
Harris~\cite{EH83}. 

\begin{theorem}\label{Th:finite}
  There are finitely many spaces of polynomials $P\in \Gr$ with a given Wronskian.
  For a general polynomial $W(t)$ of degree $(n{+}1)(d{-}n)$, there are exactly 
  $\delta(\bi_{n,d})$ spaces of polynomials with Wronskian $W(t)$.
\end{theorem}

When $W$ has distinct roots, these spaces of polynomials are exactly the points in the
intersection~\eqref{Eq:Shapiro_intersection}, where $\ba_P(s)=\bi$ at each root $s$ of
$W$.  
A limiting argument, in which the roots of the Wronskian are allowed to collide one-by-one,
proves a local form of Theorem~\ref{Th:MTV_1}. 
We say that the roots $\bs=s_1,\dotsc,s_{(n+1)(d-n)}$ of the Wronskian are 
\DeCo{{\sl clustered}} if, up to an automorphism of $\R\PP^1$, they satisfy
 \begin{equation}\label{Eq:cluster}
   0\ <\ s_1\ \ll\ s_2\ \ll\ \dotsb\ \ll\ s_{(n+1)(d-n)}\,.
 \end{equation}

\begin{theorem}[\cite{So99}]\label{Th:local}
 If the roots of a polynomial $W(t)$ of degree $(n{+}1)(d{-}n)$ are real, distinct, and
 clustered, then there are $\delta(\bi_{n,d})$ real spaces of polynomials
 with Wronskian $W(t)$ and the intersection~\eqref{Eq:Shapiro_intersection} is
 transverse. 
\end{theorem}

We noted that the intersection~\eqref{Eq:Shapiro_intersection} is not transverse when
$d=3$, $n=1$, and $W(t)=t^4-t$.
It turns out that it is always transverse when the roots of the Wronskian are distinct and
real. 
This is the stronger form of the Theorem of Mukhin, Tarasov, and Varchenko, proven
in~\cite{MTV_R}. 

\begin{theorem}\label{Th:strong}
  For any Schubert problem $\bA=(\ba^{(1)},\dotsc,\ba^{(m)})$ and any distinct real
  numbers $s_1,\dotsc,s_m$, the intersection 
 \begin{equation}\label{Eq:tr_int}
    \Omega_{\ba^{(1)}} F_\bullet(s_1)\ \bigcap\ 
    \Omega_{\ba^{(2)}} F_\bullet(s_2)\ \bigcap\ \dotsb\ \bigcap\ 
    \Omega_{\ba^{(m)}} F_\bullet(s_m)
 \end{equation}
  is transverse and consists solely of real points.
\end{theorem}

This theorem (without the transversality) is the original statement of the conjecture of
Boris Shapiro and Michael Shapiro for Grassmannians, which was posed in exactly this form to
the author in May 1995.
The Shapiro conjecture was first discussed and studied in detail in~\cite{So00}, where
significant computational evidence was presented (see also~\cite{Ve00} and~\cite{RSSS}).
These results and computations, as well as Theorem~\ref{Th:local}, highlighted the key role
that  transversality plays in the conjecture. 
Apparently, this Shapiro conjecture was in part an attempt to propose a reason for the results
in the thesis~\cite{So97} which showed that for $\G(1,d)$, there are choices of real flags
$F_\bullet^i$ in~\eqref{Eq:tr_int} so that the intersection is transverse with all points
real. 
This was extended to all problems in the special Schubert calculus on all
Grassmannians~\cite{So99}. 
Later, Vakil~\cite{Vakil} showed that this was true for all Schubert problems on all
Grassmannians. 

The  main ingredient in the proof of Theorem~\ref{Th:strong}
is an isomorphism between algebraic objects associated to the
intersection~\eqref{Eq:tr_int} and to certain representation-theoretic data.
This isomorphism provides a very deep link between Schubert calculus for the Grassmannian
and the representation theory of $\sln\C$.

We sketch the proof of Theorem~\ref{Th:MTV_1} in the next three sections.

%
\section{Spaces of polynomials with given Wronskian}\label{S:polys}

Theorem~\ref{Th:finite} enables the reduction of Theorem~\ref{Th:MTV_1} to a special
case. 
Since the Wronski map is finite, a standard limiting argument (given for example in
Section 1.3 of~\cite{MTV_Sh} or Remark 3.4 of~\cite{So00}) shows that it suffices to prove 
Theorem~\ref{Th:MTV_1} when the Wronskian has distinct real roots that are  
sufficiently general.
Since  $\delta(\bi_{n,d})$ is the upper bound for the number of spaces of polynomials 
with given Wronskian, it suffices to construct this number of distinct spaces of real
polynomials with a given Wronskian, when the Wronskian has distinct real roots that are 
sufficiently general.  
In fact, this is exactly what Mukhin, Tarasov, and Varchenko do~\cite{MTV_Sh}.\medskip 

\noindent{\bf Theorem~\ref{Th:MTV_1}\boldmath{$'$}.}
 {\it
  If $s_1,\dotsc,s_{(n{+}1)(d{-}n)}$ are generic real numbers, there are 
  $\delta(\bi_{n,d})$ distinct real vector spaces of polynomials $P$ with
  Wronskian $\prod_i(t-s_i)$.
}\medskip

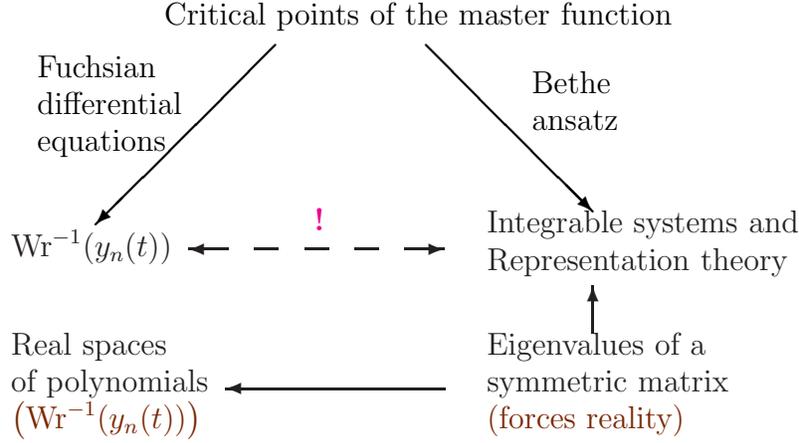
\begin{figure}[htb]
\[
  \begin{picture}(286,160)(0,15)\thicklines   

    \put(58,168){Critical points of the master function}  

    \put(100,160){\vector(-1,-1){68}}\put(157,160){\vector(1,-1){63}}

    \put(10,120){\begin{minipage}[b]{70pt}Fuchsian\newline differential\newline equations\end{minipage}}

    \put(197,128){\begin{minipage}[b]{70pt}Bethe\newline  ansatz\end{minipage}}

    \put(114.5,90){{\bf \Magenta{!}}}

    \put(  0,80){$\Wr^{-1}(y_n(t))$}

    \put( 82,83){\vector(-1,0){15}}\put( 92,83){\line(1,0){9}}
    \put(111,83){\line(1,0){9}}    \put(130,83){\line(1,0){9}}
    \put(149,83){\vector(1,0){15}}

    \put(180,75){\begin{minipage}[b]{130pt} Integrable systems and\newline
                   Representation theory\end{minipage}}

    \put(220,51){\vector(0,1){18}}

    \put(180, 15){\begin{minipage}[b]{100pt} Eigenvalues of a\newline symmetric matrix\newline
                  \Brown{(forces reality)}\end{minipage}}

    \put(164,30){\vector(-1,0){83}}

    \put( 0, 15){\begin{minipage}[b]{75pt} Real spaces\newline of polynomials\newline  
             \Brown{$\bigl(\Wr^{-1}(y_n(t))\bigr)$}\end{minipage}}

  \end{picture}
\]
\caption{Schematic of proof of Shapiro conjecture.}
\label{F:schematic}
\end{figure}

The proof first constructs $\delta(\bi_{n,d})$ distinct spaces of
polynomials with a given Wronskian having generic complex roots, 
which we describe in Section~\ref{S:construction}.
This uses a Fuchsian differential equation given by the critical points of a remarkable
symmetric function, called the master function.
The next step uses the Bethe ansatz in a certain representation $V$ of $\sln\C$: critical
points of the master function give Bethe eigenvectors of the Gaudin Hamiltonians which
turn out to be a highest weight vectors for an irreducible submodule of $V$. 
This is described in Section~\ref{S:BAGM}, where the eigenvalues of the 
Gaudin Hamiltonians on a Bethe vector are shown to be the coefficients of the Fuchsian
differential equation giving the corresponding spaces of polynomials.
This is the germ of the new, deep connection between representation theory and Schubert
calculus that led to Theorem~\ref{Th:strong}.
Finally, the Gaudin Hamiltonians are
real symmetric operators when the Wronskian has only real roots, so their
eigenvalues are real, and thus the Fuchsian differential equation
has real coefficients and the corresponding space of polynomials is also real.
Figure~\ref{F:schematic} presents a schematic of this extraordinary proof.

%
\subsection{Critical points of master functions}\label{S:construction}
The construction of $\delta(\bi_{n,d})$ spaces of polynomials with a given Wronskian
begins with the critical points of a symmetric rational function that 
arose in the study of hypergeometric solutions to the Knizhnik-Zamolodchikov
equations~\cite{SV}, and the Bethe ansatz method for the Gaudin model.

The master function depends upon parameters
$\DeCo{\bs}:=(s_1,\dotsc,s_{(n+1)(d-n)})$, which are 
the roots of our Wronskian $W$, and an additional $\binom{n+1}{2}(d{-}n)$ 
variables 
\[
    \DeCo{\bx}\ :=\ (x_1^{(1)},\dotsc,x_{d-n}^{(1)},\,
              x_1^{(2)},\dotsc,x_{2(d-n)}^{(2)},\, \dotsc\,,\,
              x_1^{(n)},\dotsc,x_{n(d-n)}^{(n)})\,.
\]
Each set of variables $\bx^{(i)}:=(x_1^{(i)},\dotsc,x_{i(d-n)}^{(i)})$ will turn out to be the
roots of certain intermediate Wronskians.

Define the \DeCo{{\sl master function} $\Phi(\bx;\bs)$} by the (rather formidable) formula
 \begin{equation}\label{Eq:MasterFunction}
   \frac{\displaystyle \prod_{i=1}^n\ 
    \prod_{1\leq j<k\leq i(d-n)}(x_j^{(i)}-x_k^{(i)})^2
    \ \cdot\ \prod_{1\leq j<k<(n+1)(d-n)}(s_j-s_k)^2}{\displaystyle
   \prod_{i=1}^{n-1}\ \prod_{j=1}^{i(d-n)}\ \prod_{k=1}^{(i+1)(d-n)}
   (x_j^{(i)}-x_k^{(i+1)})\ \cdot\ 
   \prod_{j=1}^{n(d-n)}\ \prod_{k=1}^{(n+1)(d-n)}(x_j^{(n)}-s_k)}\ .
 \end{equation}
This is separately symmetric in each set of variables $\bx^{(i)}$.
The  Cartan matrix for $\sln$ appears in the exponents of the factors
$(x_*^{(i)}-x_*^{(j)})$ in~\eqref{Eq:MasterFunction}.
This hints at the relation of these master functions to Lie theory, which
we do not discuss.

The critical points of the master function are solutions to the system of equations
 \begin{equation}\label{Eq:Critical}
   \frac{1}{\Phi}\frac{\partial}{\partial x_j^{(i)}}\Phi(\bx;\bs)\ =\ 0
   \qquad\mbox{for}\quad i=1,\dotsc,n,\quad
                         j=1,\dotsc,i(d{-}n)\,.
 \end{equation}
When the parameters $\bs$ are generic, these \DeCo{{\sl Bethe ansatz equations}} turn out
to have finitely many solutions. 
The master function is invariant under the group 
\[
   \DeCo{\calS}\ :=\ 
    \calS_{d-n}\times\calS_{2(d-n)}\times\,\dotsb\,\times\calS_{n(d-n)}\,,
\]
where $\calS_m$ is the group of permutations of $\{1,\dotsc,m\}$, and 
the factor $\calS_{i(d{-}n)}$ permutes the variables in $\bx^{(i)}$.
Thus $\calS$ acts on the critical points.
The invariants of this action are polynomials whose roots are 
the coordinates of the critical points.

Given a critical point $\bx$, define monic polynomials $\bp_\bx:=(p_1,\dotsc,p_n)$ 
where the components $\bx^{(i)}$ of $\bx$ are the roots of $p_i$,
 \begin{equation}\label{eq:crit_polys}
    \DeCo{p_i}\ :=\ \prod_{j=1}^{i(d-n)} (t-x_j^{(i)})
    \qquad\mbox{for}\quad i=1,\dotsc,n\,.
 \end{equation}
Also write $\DeCo{p_{n+1}}$ for the Wronskian, the monic polynomial with roots $\bs$.
The discriminant \DeCo{$\mbox{Discr}(f)$} of a polynomial $f$ is the square of
the product of differences of its roots and the resultant \DeCo{$\mbox{Res}(f,g)$} is the
product of all differences of the roots of $f$ and $g$~\cite{CLO}.
Then the formula for the master function~\eqref{Eq:MasterFunction} becomes
 \begin{equation}\label{Eq:Resultant}
  \Phi(\bx;\bs)\ =\ \prod_{i=1}^{n+1} \mbox{Discr}(p_i) 
                    \Bigg/  \prod_{i=1}^n \mbox{Res}(p_i,p_{i+1}) \ .
 \end{equation}

The connection between the critical points of $\Phi(\bx;\bs)$ and spaces of
polynomials with Wronskian $W$ is through a Fuchsian differential equation.
Given (an orbit of) a critical point $\bx$ represented by the list of polynomials
$\bp_\bx$, define the
 \DeCo{{\sl fundamental differential operator $D_\bx$ of the critical point $\bx$}} by
 \begin{equation}\label{Eq:FundDiffOp}
   \Bigl(\frac{d}{dt} - \dlog\Bigl(\frac{W}{p_n}\Bigr)\Bigr)
   \,\dotsb\,
   \Bigl(\frac{d}{dt} - \dlog\Bigl(\frac{p_{2}}{p_1}\Bigr)\Bigr)
   \Bigl(\frac{d}{dt} - \dlog(p_1)\Bigr)\,,
 \end{equation}
where $\dlog(f):=\frac{d}{dt}\ln f$.
The kernel \DeCo{$V_\bx$} of $D_\bx$ is the 
\DeCo{{\sl fundamental space of the critical point $\bx$}}.

\begin{example}\label{Ex:solutions}
 Since 
\[
   \Bigr( \frac{d}{dt}\ -\ \dlog(p)\Bigr)\; p\ =\ 
   \Bigr( \frac{d}{dt}\ -\ \frac{p'}{p}\Bigr)\; p\ =\ 
   p' - \frac{p'}{p}p\ =\ 0\,,
\]
 we see that $p_1$ is a solution of $D_\bx$.
 It is instructive to look at $D_\bx$ and $V_\bx$ when $n=1$.
 Suppose that $f$ a solution to $D_\bx$ that is linearly independent from $p_1$.
 Then
\[
  0\ =\ 
   \Bigl(\frac{d}{dt} - \dlog\Bigl(\frac{W}{p_{1}}\Bigr)\Bigr)
   \Bigl(\frac{d}{dt} - \dlog(p_1)\Bigr)\;f
  \ =\  \Bigl(\frac{d}{dt} - \dlog\Bigl(\frac{W}{p_{1}}\Bigr)\Bigr)
     \bigl( f' - \frac{p_1'}{p_1} f \bigr)\,.
\]
 This implies that
\[
   \frac{W}{p_1}\ =\ f' - \frac{p_1'}{p_1} f\,,
\]
 so $W=\Wr(f,p_1)$, and the kernel of $D_\bx$ is a 2-dimensional space of functions with
 Wronskian $W$. 
\end{example}

What we just saw is always the case.
The following result is due to Scherbak and Varchenko~\cite{ScVa} for $n=1$ and to Mukhin
and Varchenko~\cite[\S 5]{MV04} 
for all $n$.

\begin{theorem}\label{Th:fundSpace}
 Suppose that $V_\bx$ is the fundamental space of a critical point $\bx$ of
 the master function $\Phi$ with generic parameters $\bs$ which are the roots of $W$. 
 \begin{enumerate}
  \item Then $V_\bx$ is an $(n{+}1)$-dimensional space of polynomials of degree $d$ lying in 
          $\Gr^\circ$ with Wronskian $W$.
  \item The critical point $\bx$ is recovered from $V_\bx$ in some cases as follows.
        Suppose that $f_0,\dotsc,f_n$ are monic polynomials in $V_\bx$ with 
        $\deg f_i=d{-}n+i$, each $f_i$ is square-free, and that the pairs $f_i$ and
        $f_{i+1}$ are relatively prime.
        Then, up to scalar multiples, the polynomials $p_1,\dotsc,p_n$ in the sequence
        $\bp_\bx$ are
\[
    f_0\,,\ \Wr(f_0,f_1)\,,\ \Wr(f_0,f_1,f_2)\,,\ \dotsc\,,\ 
    \Wr(f_0,\dotsc,f_{n-1})\,.
\]
 \end{enumerate} 
\end{theorem}

Statement~(2) includes a general result about factoring a linear differential operator into
differential operators of degree 1.
Linearly independent  $C^\infty$ functions $f_0,\dotsc,f_n$ 
span the kernel of the differential operator of degree $n{+}1$,
\[
   \det\left(\begin{matrix}
          f_0 & f_1 & \dotsb& f_n &  1 \\
          f'_0& f'_1& \dotsb& f'_n&\frac{d}{dt}\\
        \vdots&\vdots&\ddots&\vdots&\vdots\\
     f_0^{(n+1)}& f_1^{(n+1)}& \dotsb& f_n^{(n+1)}&\frac{d^{n+1}}{dt^{n+1}}
    \end{matrix}\right)\ .
\]
If we set $p_{i+1}:=\Wr(f_0,\dotsc,f_i)$, then~\eqref{Eq:FundDiffOp} is a factorization
over $\C(t)$ of this determinant into differential operators of degree 1.
This follows from some interesting identities among Wronskians shown in the Appendix
of~\cite{MV04}.

Theorem~\ref{Th:fundSpace} is deeper than this curious fact.
When the polynomials $p_1,\dotsc,p_n,W$ are square-free,
consecutive pairs are relatively prime, and $\bs$ is generic, it implies that the kernel $V$ of
an operator of the form~\eqref{Eq:FundDiffOp} is a space of polynomials with Wronskian $W$
having roots $\bs$ if and only if the polynomials $p_1,\dotsc,p_n$ come from the
critical points of the master function~\eqref{Eq:MasterFunction} corresponding to $W$.

This gives an injection from $\calS$-orbits of critical
points of the master function $\Phi$ with parameters $\bs$ to spaces of polynomials in
$\Gr^\circ$ whose Wronskian has roots $\bs$.
Mukhin and Varchenko showed that this is a bijection when $\bs$ is generic. 

\begin{theorem}[Theorem 6.1 in~\cite{MV05}]\label{Th:distinct}
  For generic complex numbers $\bs$, the master function $\Phi$ has nondegenerate critical
  points that form $\delta(\bi_{n,d})$ distinct orbits. 
\end{theorem}

The structure (but not of course the details) of their proof is remarkably similar to the
structure of the proof of Theorem~\ref{Th:local}; they allow the
parameters to collide one-by-one, and study how the orbits of critical points behave.
Ultimately, they obtain the same recursion as in~\cite{So99}, which mimics the Pieri
formula for the branching rule for tensor products of representations of $\sln$ with its
fundamental representation $V_{\omega_n}$.
This same structure is also found in the main argument in~\cite{EG02c}.
In fact, this is the same recursion in $\ba$ that Schubert established for intersection
numbers $\delta(\ba,\bi,\dotsc,\bi)$, and then solved to obtain the
formula~\eqref{Eq:Grass_degree} in~\cite{Sch1886c}. 

%
\section{The Bethe  ansatz for the Gaudin model}\label{S:BAGM}

The Bethe ansatz is a general (conjectural) method to find pure states, called 
\DeCo{{\sl Bethe vectors}}, of quantum integrable systems. 
The (periodic) Gaudin model is an integrable system consisting of a family of commuting
operators called the Gaudin Hamiltonians that act on a representation $V$ of $\sln\C$.
In this Bethe ansatz, a vector-valued rational function 
is constructed so that for certain values of the parameters it yields a complete
set of Bethe vectors.
As the Gaudin Hamiltonians commute with the action of $\sln\C$, the Bethe vectors turn out to
be highest weight vectors generating irreducible submodules of $V$, and so this also gives a
method for decomposing some representations $V$ of $\sln\C$ into irreducible submodules. 
The development, justification, and refinements of this Bethe ansatz are the subject of a large
body of  work, a small part of which we mention.

%
\subsection{Representations of  $\mathfrak{sl}_{n+1}\C$}
The Lie algebra $\sln\C$ (or simply $\sln$) is the space of $(n{+}1)\times(n{+}1)$-matrices
with trace zero.
It has a decomposition
\[
   \sln\ =\ \frakn_-\oplus\frakh\oplus\frakn_+\,,
\]
where $\frakn_+$ $(\frakn_-)$ are the strictly upper (lower) triangular matrices,
and $\frakh$ consists of the diagonal matrices with zero trace.
The universal enveloping algebra $U\sln$ of $\sln$ is the associative algebra generated by
$\sln$ subject to the relations $uv-vu=[u,v]$ for $u,v\in\sln$ where $[u,v]$ is the Lie
bracket in $\sln$.

We consider only finite-dimensional representations of $\sln$ (equivalently, of $U\sln$).
For a more complete treatment, see~\cite{FuHa}.
Any representation $V$ of $\sln$ decomposes into joint
eigen\-spaces of $\frakh$, called \DeCo{{\sl weight spaces}},
\[
   V\ =\ \bigoplus_{\mu\in\frakh^*} V[\mu]\,,
\]
where, for $v\in V[\mu]$ and $h\in\frakh$, we have $h.v=\mu(h)v$.
The possible weights $\mu$ of representations lie in the integral 
\DeCo{{\sl weight lattice}}.
This has a distinguished basis of \DeCo{{\sl fundamental weights}}
$\DeCo{\omega_1},\dotsc,\DeCo{\omega_n}$ that generate the cone of 
\DeCo{{\sl dominant weights}}.

An irreducible representation $V$ has a unique one-dimensional weight space that is annihilated
by the nilpotent subalgebra $\frakn_+$ of $\sln$.
The associated weight $\mu$ is dominant, and it is called the \DeCo{{\sl highest weight}} of $V$.
Any nonzero vector with this weight is a highest weight vector of $V$, and it generates $V$.
Furthermore, any two irreducible modules with the same highest weight are isomorphic.
Write \DeCo{$V_\mu$} for the \DeCo{{\sl highest weight module}} with highest weight $\mu$.
Lastly, there is one highest weight module for each dominant weight.

More generally, if $V$ is any representation of $\sln$ and $\mu$ is a weight, 
then the \DeCo{{\sl singular vectors}} in $V$ of weight $\mu$, 
written \DeCo{$\sing(V[\mu])$}, are the vectors in $V[\mu]$ annihilated by $\frakn_+$. 
If $v\in \sing(V[\mu])$ is nonzero, then the submodule $U\sln.v$ it generates is
isomorphic to the highest weight module $V_\mu$.
Thus $V$ decomposes as a direct sum of submodules generated by the singular
vectors, 
 \begin{equation}\label{Eq:sing_decomp}
   V\ =\ \bigoplus_\mu  U\sln.\sing (V[\mu])\,,
 \end{equation}
so that the multiplicity of the highest weight module $V_\mu$ in $V$ is simply the
dimension of its space of singular vectors of weight $\mu$.

When $V$ is a tensor product of highest weight modules, the Littlewood-Richard\-son
rule~\cite{Fu} gives formulas for the dimensions of the spaces of singular vectors.
Since this is the same rule for the number of points in an
intersection~\eqref{Eq:Schubert_intersection} of Schubert varieties from a
Schubert problem,  these geometric intersection numbers are equal to the dimensions of 
spaces of singular vectors.
In particular, if $V_{\omega_1}\simeq\C^{n+1}$ is the defining representation of $\sln$
and $\DeCo{V_{\omega_n}}=\wedge^n V_{\omega_1}=V_{\omega_1}^*$ (these are the first and last
fundamental representations of $\sln$), then 
 \begin{equation}\label{Eq:SC=RT}
   \dim \sing( V_{\omega_n}^{\otimes (n+1)(d-n)}[0])\ =\ 
    \delta(\bi_{n,d})\,.
 \end{equation}
It is important to note that this equality of numbers is purely formal, in that the same
formula governs both numbers.
A direct connection remains to be found.

%
\subsection{The (periodic) Gaudin model}
The Bethe ansatz is a conjectural method to obtain a complete set of eigenvectors for the
integrable system on $V:=V_{\omega_n}^{\otimes m}$ given by the Gaudin Hamiltonians (defined
below). 
Since these Gaudin Hamiltonians commute with $\sln$, the Bethe ansatz has the additional
benefit of giving an explicit basis for $\sing(V[\mu])$, thus explicitly giving the 
decomposition~\eqref{Eq:sing_decomp}.

The Gaudin Hamiltonians act on $V_{\omega_n}^{\otimes m}$ and depend upon $m$
distinct complex numbers $s_1,\dotsc,s_m$ and a complex variable $t$.
Let $\gln$ be the Lie algebra of $(n{+}1)\times(n{+}1)$ complex matrices.
For each $i,j=1,\dotsc,n{+}1$, let $E_{i,j}\in\gln$ be the matrix whose only nonzero entry is a
1 in row $i$ and column $j$.
For each pair $(i,j)$ consider the differential operator $X_{i,j}(t)$ acting on 
$V_{\omega_n}^{\otimes m}$-valued functions of $t$,
\[
   \DeCo{X_{i,j}(t)}\ :=\ \delta_{i,j}\frac{d}{dt}\ -\ 
    \sum_{k=1}^m \frac{E_{j,i}^{(k)}}{t-s_k}\ ,
\]
where $E_{j,i}^{(k)}$ acts on tensors in $V_{\omega_n}^{\otimes m}$  by $E_{j,i}$ in the
$k$th factor and by the identity in other factors.
Define a differential operator acting on $V_{\omega_n}^{\otimes m}$-valued functions of
$t$, 
\[
   \DeCo{\bM}\ := \sum_{\sigma\in\calS} \mbox{sgn}(\sigma)\; 
    X_{1,\sigma(1)}(t)\; 
    X_{2,\sigma(2)}(t)\; \dotsb\; X_{n+1,\sigma(n+1)}(t)\ ,
\]
where $\calS$ is the group of permutations of $\{1,\dotsc,n{+}1\}$ and
$\mbox{sgn}(\sigma)=\pm$ is the sign of a permutation $\sigma\in\calS$.
Write $\bM$ in standard form
\[
   \bM\ =\ \frac{d^{n+1}}{dt^{n+1}}\ +\ M_1(t) \frac{d^n}{dt^n}\ +\ 
    \dotsb\ +\ M_{n+1}(t)\,.
\]
These coefficients $M_1(t),\dotsc,M_{n+1}(t)$ are called the (higher) 
\DeCo{{\sl Gaudin Hamiltonians}}.
They are linear operators that depend rationally on $t$ and act on  
$V_{\omega_n}^{\otimes m}$.
We collect together some of their properties.

\begin{theorem}
 Suppose that $s_1,\dotsc,s_m$ are distinct complex numbers.
 Then
 \begin{enumerate}
   \item The Gaudin Hamiltonians commute, that is, $[M_i(u),M_j(v)]=0$ for all
      $i,j=1,\dotsc,n{+}1$ and $u,v\in\C$.
   \item The Gaudin Hamiltonians commute with the action of\/ $\sln$ on
      $V_{\omega_n}^{\otimes m}$.
 \end{enumerate}
\end{theorem}

Proofs are given in~\cite{KuS}, as well as
Propositions 7.2 and 8.3 in~\cite{MTV_06}, and are based on results of Talalaev~\cite{Ta}.
A consequence of the second assertion is that the Gaudin Hamiltonians preserve the weight
space decomposition of the singular vectors of $V_{\omega_n}^{\otimes m}$.
Since they commute, the singular vectors of $V_{\omega_n}^{\otimes m}$
have a basis of common eigenvectors of the Gaudin Hamiltonians.
The Bethe ansatz is a method to write down joint eigenvectors and their eigenvalues. 

%
%
\subsection{The Bethe ansatz for the Gaudin model}\label{S:BAGMwv}

This begins with a rational function that takes
values in a weight space $V_{\omega_n}^{\otimes m}[\mu]$, 
\[
   v\ \colon\ \C^l\times\C^m\ \longmapsto\  V_{\omega_n}^{\otimes m}[\mu]\ .
\]
This \DeCo{{\sl universal weight function}} was introduced in~\cite{SV} to solve the
Knizhnik-Zamo\-lod\-chi\-kov equations with values in $V_{\omega_n}^{\otimes m}[\mu]$.
When $(\bx,\bs)$ is a critical point of a master function, the vector 
$v(\bx,\bs)$ is both singular and an eigenvector of the Gaudin Hamiltonians.
(This master function is a generalization of the one defined
by~\eqref{Eq:MasterFunction}.) 
The Bethe ansatz conjecture for the periodic Gaudin model asserts that the 
vectors $v(\bx,\bs)$ form a basis for the space of singular vectors.

Fix a highest weight vector $\DeCo{v_{n+1}}\in V_{\omega_n}[\omega_n]$.
Then $v_{n+1}^{\otimes m}$ generates $V_{\omega_n}^{\otimes m}$ as a
$U\sln^{\otimes m}$-module.  
In particular, any vector in $V_{\omega_n}^{\otimes m}$ is a linear combination of vectors that
are obtained from $v_{n+1}^{\otimes m}$ by applying a sequence of operators $E^{(k)}_{i+1,i}$, for
$1\leq k\leq m$ and $1\leq i\leq n$.
The universal weight function is a linear combination of such vectors of weight
$\mu$.

When 
$m=(n{+}1)(d{-}n)$, $l=\binom{n+1}{2}(d{-}n)$, and $\mu=0$,
the universal weight function is a map
\[
  v\ \colon\ \C^{\binom{n+1}{2}(d{-}n)}\times\C^{(n+1)(d{-}n)} \ 
    \longrightarrow\  V_{\omega_n}^{\otimes(n+1)(d-n)}[0]\,.
\]
To describe it, note that a vector $E_{a+1,a}E_{b+1,b}\dotsb E_{c+1,c}.v_{n+1}$ is nonzero only if
\[
  (a,b,\dotsc,c)\ =\ (a, a{+}1,\dotsc,n{-}1,n)\,.
\]
Write \DeCo{$v_a$} for this vector.
The vectors $v_1,\dotsc,v_{n{+}1}$ form a basis of $V_{\omega_n}$.
Thus only some sequences of operators $E^{(k)}_{i+1,i}$ applied to $v_{n+1}^{\otimes(n+1)(d-n)}$ 
give a nonzero vector.
These sequences are completely determined once we know the weight of the result.
The operator $E^{(k)}_{i+1,i}$ lowers the weight of a weight vector by the root $\alpha_i$.
Since
 \begin{equation}\label{eq:weight_sum}
    (n{+}1)\omega_n\ =\ \alpha_1+2\alpha_2+\dotsb+n\alpha_n\,,
 \end{equation}
there are $i(d{-}n)$ occurrences of $E^{(k)}_{i+1,i}$, which is the number of
variables in $\bx^{(i)}$.

Let \DeCo{$\calB$} be the set of all sequences $(b_1,b_2,\dotsc,b_{(n+1)(d-n)})$,
where $1\leq b_k\leq n{+}1$ for each $k$ and we have
\[
   \#\{k\mid b_k\leq i\}\ =\ i(d{-}n)\,.
\]
Given a sequence $B$ in $\calB$, define
 \begin{eqnarray*}
   \DeCo{v_B} &:=& v_{b_1}\otimes v_{b_2}\otimes\dotsb\otimes v_{b_{(n+1)(d-n)}}\\
               &=& \bigotimes_{k=1}^{(n+1)(d-n)} 
           \bigl( E^{(k)}_{b_k+1,b_k} \dotsb E^{(k)}_{n,n-1}\cdot
           E^{(k)}_{n+1,n}\bigr).v_{n+1}\,,
 \end{eqnarray*}
where the operator $ E^{(k)}_{b_k+1,b_k}\dotsb E^{(k)}_{n,n-1}\cdot E^{(k)}_{n+1,n}$ is
the identity if $b_k=n+1$.
Then $v_B$ is a vector of weight $0$, by~\eqref{eq:weight_sum}.
The universal weight function is a linear combination of these vectors $v_B$,
\[
   v(\bx;\bs)\ =\ \sum_{B\in\calB} w_B(\bx;\bs)\cdot v_B\,,
\]
where $w_B(\bx,\bs)$ is separately symmetric in each set of variables
$\bx^{(i)}$. 

To describe $w_B(\bx;\bs)$, suppose that 
\[
   \bz\ =\ (\bz^{(1)},\bz^{(2)},\dotsc,\bz^{((n+1)(d-n))})
\]
is a partition of the variables $\bx$ into $(n{+}1)(d{-}n)$ sets of variables where the $k$th
set $\bz^{(k)}$ of variables has exactly one variable from each set $\bx^{(i)}$ with $b_k\leq i$
(and is empty when $b_k=n{+}1$).
That is, if $b_k\leq n$, then
 \begin{equation}\label{eq:yk}
   \bz^{(k)}\ =\ (x^{(b_k)}_{c_{b_k}}, x^{(b_k+1)}_{c_{b_k+1}},\dotsc, x^{(n)}_{c_n})\,,
 \end{equation}
for some indices $c_{b_k},\dotsc,c_n$. 
If $b_k=n{+}1$, set $w^{(k)}(\bz):=1$, and otherwise
\[
   w^{(k)}(\bz;\bs)\ :=\ \frac{1}{x^{(b_k)}_{c_{b_k}}-x^{(b_k+1)}_{c_{b_k+1}}}\dotsb
               \frac{1}{x^{(n-1)}_{c_{n-1}}-x^{(n)}_{c_{n}}}\cdot
               \frac{1}{x^{(n)}_{c_{n}}-s_k}\,,
\]
in the notation~\eqref{eq:yk}.
Then we set 
\[
   w(\bz;\bs)\ :=\ \prod_{k=1}^{(n+1)(d-n)} w^{(k)}(\bz;\bs)\,.
\]
Finally, $w_B(\bx;\bs)$ is the sum of the rational functions $w(\bz;\bs)$ over all 
such partitions $\bz$ of the variables $\bx$.
(Equivalently, the symmetrization of any single $w(\bz;\bs)$.)

While $v(\bx,\bs)$ is a rational function of $\bx$ and hence not globally defined, 
if the coordinates of $\bs$ are distinct and $\bx$ is a critical point of the master
function~\eqref{Eq:MasterFunction}, then the vector 
$v(\bx,\bs)\in V^{\otimes(n{+}1)(d{-}n)}_{\omega_n}[0]$ is well-defined, nonzero and it is in fact
a singular vector (Lemma~2.1 of~\cite{MV05}).
Such a vector $v(\bx,\bs)$ when $\bx$ is a critical point of the master function is called a 
\DeCo{{\sl Bethe vector}}.
Mukhin and Varchenko also prove the following, which is the second part of Theorem~6.1
in~\cite{MV05}.

\begin{theorem}\label{Th:BV_basis}
 When $\bs\in\C^{(n{+}1)(d-n)}$ is general, the Bethe vectors form a basis of  the
 space $\sing\bigl( V^{\otimes(n{+}1)(d{-}n)}_{\omega_n}[0]\bigr)$.
\end{theorem}

These Bethe vectors are the joint
eigenvectors of the Gaudin Hamiltonians.

\begin{theorem}[Theorem 9.2 in \cite{MTV_06}]\label{Th:BV_EV}
 For any critical point $\bx$ of the master function~$\eqref{Eq:MasterFunction}$, the
 Bethe vector $v(\bx,\bs)$ is a joint eigenvector of the Gaudin Hamiltonians
 $M_1(t),\dotsc,M_{n+1}(t)$.
 Its eigenvalues $\mu_1(t),\dotsc,\mu_{n+1}(t)$ are given by the formula
 \begin{multline}\label{Eq:ev}
   \quad \frac{d^{n+1}}{dt^{n+1}}\ +\ \mu_1(t)\frac{d^n}{dt^n}\ +\ 
      \dotsb\ +\ \mu_n(t)\frac{d}{dt}\ +\ \mu_{n+1}(t)\ =\ \\ \rule{0pt}{17pt}
     \Bigl(\frac{d}{dt}+\dlog(p_1)\Bigr)
     \Bigl(\frac{d}{dt}+\dlog\Bigl(\frac{p_2}{p_1}\Bigr)\Bigr)
       \ \dotsb\ 
     \Bigl(\frac{d}{dt}+\dlog\Bigl(\frac{p_n}{p_{n-1}}\Bigr)\Bigr)
     \Bigl(\frac{d}{dt}+\dlog\Bigl(\frac{W}{p_n}\Bigr)\Bigr)\,,\quad
 \end{multline}
 where $p_1(t),\dotsc,p_n(t)$ are the polynomials~\eqref{eq:crit_polys} associated to
 the critical point $\bx$ and $W(t)$ is the polynomial with roots $\bs$.
\end{theorem}

Observe that~\eqref{Eq:ev} is similar to the formula~\eqref{Eq:FundDiffOp} for the
differential operator $D_\bx$ of the critical point $\bx$.
This similarity is made more precise if we replace the Gaudin Hamiltonians by a different
set of operators.
Consider the differential operator formally conjugate to $(-1)^{n+1}M$,
 \begin{eqnarray*}
   \DeCo{K}&=& \frac{d^{n+1}}{dt^{n+1}}\ -\ \frac{d^n}{dt^n}M_1(t)\ +\ 
    \dotsb\ +\ (-1)^n\frac{d}{dt} M_n(t)\ +\ (-1)^{n+1} M_{n+1}(t)\\
   &=& \frac{d^{n+1}}{dt^{n+1}}\ +\ \DeCo{K_1(t)}\frac{d^n}{dt^n}\ +\ 
    \dotsb\ +\  \DeCo{K_n(t)}\frac{d}{dt}\ +\ \DeCo{K_{n+1}(t)}\ .
 \end{eqnarray*}
These coefficients $K_i(t)$ are operators on $V_{\omega_n}^{\otimes(n+1)(d-n)}$
that depend rationally on $t$, and are also called the Gaudin Hamiltonians.
Here are the first three, 
 \begin{eqnarray*}
   K_1(t)&=&-M_1(t)\,,\qquad\qquad\qquad
   K_2(t)\ =\ M_2(t)\ -\ n M_1'(t)\,,\\
   K_3(t)&=& -M_3(t)\ +\ (n{-}1)M_2''(t)\ -\ \binom{n}{2}M_1'''(t)\,,
 \end{eqnarray*}
and in general $K_i(t)$ is a differential polynomial in $M_1(t),\dotsc,M_i(t)$.

Like the $M_i(t)$, these operators commute with each other
and with $\sln$, and the Bethe vector $v(\bx,\bs)$ is a joint eigenvector of these new Gaudin 
Hamiltonians $K_i(t)$. 
The corresponding eigenvalues $\lambda_1(t),\dotsc,\lambda_{n{+}1}(t)$ are given
by the formula
 \begin{multline}\label{Eq:Kev}
   \quad \frac{d^{n+1}}{dt^{n+1}}\ +\ \lambda_1(t)\frac{d^n}{dt^n}\ +\ 
      \dotsb\ +\ \lambda_n(t)\frac{d}{dt}\ +\ \lambda_{n+1}(t)\ =\ \\
   \Bigl(\frac{d}{dt} - \dlog\Bigl(\frac{W}{p_{n}}\Bigr)\Bigr)   \rule{0pt}{17pt}
   \Bigl(\frac{d}{dt} - \dlog\Bigl(\frac{p_n}{p_{n-1}}\Bigr)\Bigr)
   \,\dotsb\,
   \Bigl(\frac{d}{dt} - \dlog\Bigl(\frac{p_{2}}{p_1}\Bigr)\Bigr)
   \Bigl(\frac{d}{dt} - \dlog(p_0)\Bigr)\,,
 \end{multline}
which is $(!)$ the fundamental differential operator $D_\bx$ of the critical point $\bx$.  

\begin{corollary}\label{Co:GH_simple}
 Suppose that $\bs\in\C^{(n+1)(d-n)}$ is generic.
 \begin{enumerate}
  \item  The Bethe vectors form an eigenbasis of\/ 
       $\sing(V^{\otimes(n+1)(d-n)}_{\omega_n}[0])$ for the Gaudin Hamiltonians 
       $K_1(t),\dotsc,K_{n+1}(t)$.
  \item  The Gaudin Hamiltonians $K_1(t),\dotsc,K_{n+1}(t)$ have simple spectrum in that
       their eigenvalues separate the basis of eigenvectors.
\end{enumerate}
\end{corollary}

Statement (1) follows from Theorems~\ref{Th:BV_basis} and~\ref{Th:BV_EV}.
For Statement (2), suppose that two Bethe vectors $v(\bx,\bs)$ and $v(\bx',\bs)$ have the
same eigenvalues.
By~\eqref{Eq:Kev}, the corresponding fundamental differential operators would be equal,
$D_{\bx}=D_{\bx'}$.
But this implies that the fundamental spaces coincide, $V_\bx=V_{\bx'}$.
By Theorem~\ref{Th:fundSpace} the fundamental space determines the orbit of critical
points, so the critical points $\bx$ and $\bx'$ lie in the same orbit, which implies that 
 $v(\bx,\bs)=v(\bx',\bs)$.

%
%
\section{Shapovalov form and the proof of the Shapiro conjecture}\label{S:Shapovalov}
The last step in the proof of Theorem~\ref{Th:MTV_1} is to show that if
$\bs\in\R^{(n+1)(d-n)}$ is generic and $\bx$ is a critical point of the master
function~\eqref{Eq:MasterFunction}, then the fundamental space $V_\bx$ of the critical
point $\bx$ has a basis of real polynomials.
The reason for this reality is that the eigenvectors and
eigenvalues of a symmetric matrix are real.

We begin with the Shapovalov form.
The map $\tau\colon E_{i,j}\mapsto E_{j,i}$ induces an antiautomorphism on $\sln$.
Given a highest weight module $V_\mu$, and a highest weight vector $v\in V_\mu[\mu]$,
the \DeCo{{\sl Shapovalov form} $\langle \cdot,\cdot\rangle$} on $V_\mu$ is defined
recursively by 
\[
  \langle v,v\rangle\ =\ 1\qquad\mbox{and}\qquad
  \langle g.u,v\rangle\ =\ \langle u,\tau(g).v\rangle\,,
\]
for $g\in\sln$ and $u,v\in V$.
In general, this Shapovalov form is nondegenerate on $V_\mu$ and positive definite on the
real part of $V_\mu$.

For example, the Shapovalov form on $V_{\omega_n}$ is the standard Euclidean inner
product, $\langle v_i,v_j\rangle=\delta_{i,j}$, in the basis $v_1,\dotsc,v_{n{+}1}$ of
Section~\ref{S:BAGMwv}.
This induces the symmetric (tensor)
Shapovalov form on the tensor product $V_{\omega_n}^{\otimes(n+1)(d-n)}$,
which is  positive definite on the real part of 
$V_{\omega_n}^{\otimes(n+1)(d-n)}$.

\begin{theorem}[Proposition 9.1 in~\cite{MTV_06}]
  The Gaudin Hamiltonians are symmetric with respect to the tensor Shapovalov form,
\[
   \langle K_i(t).u,\, v\rangle \ =\ \langle u,\, K_i(t).v\rangle\,,
\]
  for all $i=1,\dotsc,n{+}1$, $t\in\C$, and $u,v\in V_{\omega_n}^{\otimes(n+1)(d-n)}$.
\end{theorem}

We give the most important consequence of this result for our story.

\begin{corollary}\label{C:symmetric}
  When the parameters $\bs$ and variable $t$ are real, the Gaudin Hamiltonians 
  $K_1(t),\dotsc, K_{n+1}(t)$ are real linear operators with real spectrum.
\end{corollary}

\noindent{\it Proof.}
 The Gaudin Hamiltonians $M_1(t),\dotsc,M_{n+1}(t)$ are real linear operators which
 act on the real part of  $V_{\omega_n}^{\otimes(n+1)(d-n)}$, by their definition.
 The same is then also true of the Gaudin Hamiltonians $K_1(t),\dotsc,K_{n+1}(t)$.
 But these are symmetric with respect to the Shapovalov form
 and thus have real spectrum.
\hfill\raisebox{-4.5pt}{\includegraphics[height=14pt]{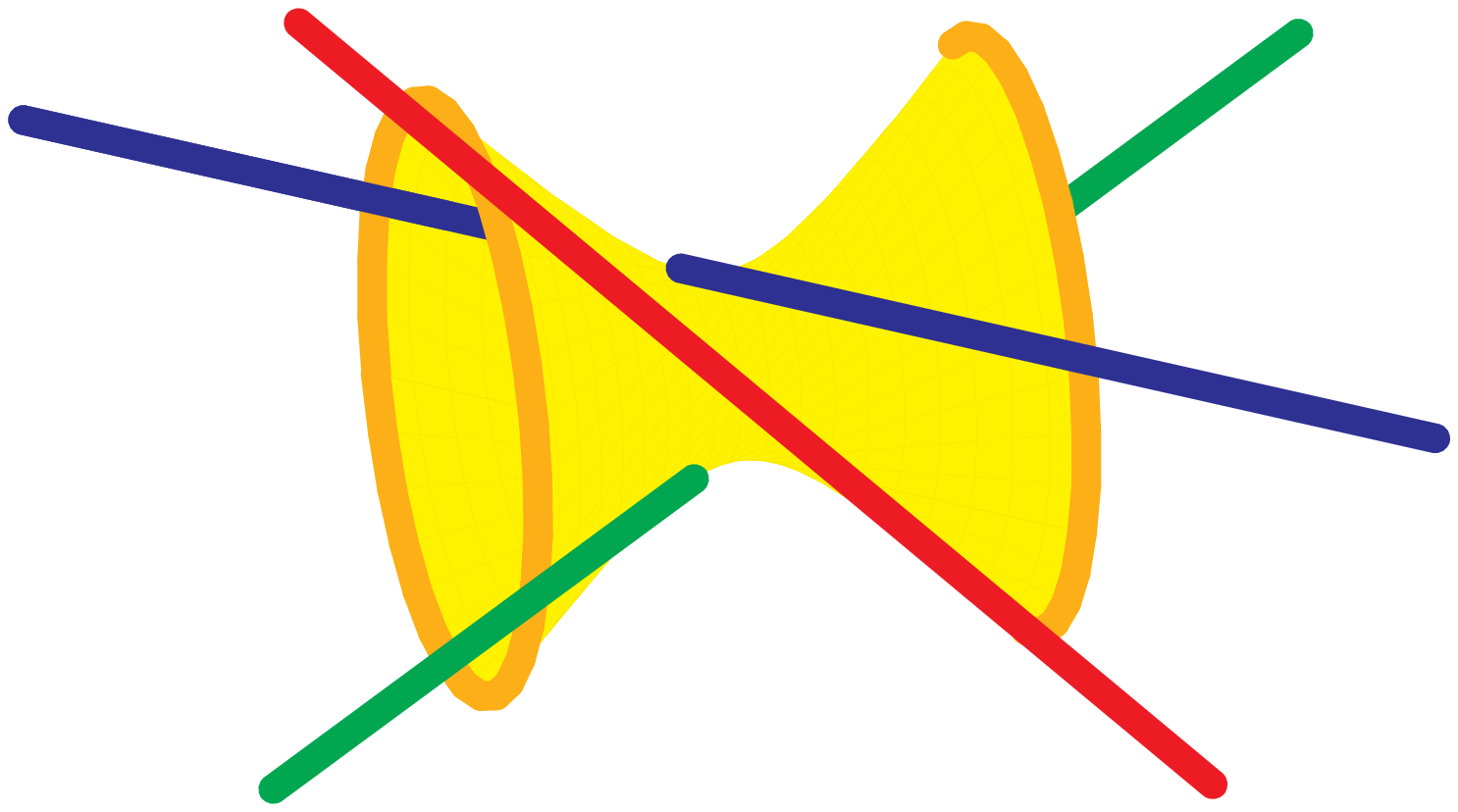}}\vspace{5pt}

\noindent{\it Proof of Theorem~$\ref{Th:MTV_1}$.}
 Suppose that $\bs\in\R^{(n+1)(d-n)}$ is general.
 By Corollary~\ref{C:symmetric}, the Gaudin Hamiltonians for $t\in\R$ acting on 
 $\sing(V_{\omega_n}^{(n+1)(d-n)}[0])$ are symmetric operators on a Euclidean space, and
 so have real eigenvalues.
 The Bethe vectors $v(\bx,\bs)$ for critical points $\bx$ of the master function with
 parameters $\bs$ form an eigenbasis for the Gaudin Hamiltonians.
 As $\bs$ is general, the eigenvalues are distinct by Corollary~\ref{Co:GH_simple} (2),
 and so the Bethe vectors must be real.

 Given a critical point $\bx$, the eigenvalues $\lambda_1(t),\dotsc,\lambda_{n+1}(t)$ of
 the Bethe vectors are then real rational functions, and so the fundamental differential
 operator $D_\bx$ has real coefficients.
 But then the fundamental space $V_\bx$ of polynomials is real.

 Thus each of the $\delta(\bi_{n,d})$ spaces of polynomials $V_\bx$
 whose Wronskian has roots $\bs$ that were constructed in Section~\ref{S:polys} is in fact
 real.
 This proves Theorem~\ref{Th:MTV_1}.
\hfill\raisebox{-4.5pt}{\includegraphics[height=14pt]{pictures/qed.eps}}\smallskip


%
%
\section{Other proofs of the Shapiro conjecture}\label{S:other}
The proofs of different Bethe ans\"atze for other models (other integrable systems) and
other Lie algebras, which is ongoing work of Mukhin, Tarasov, and Varchenko, and others,
can lead to generalizations of Theorem~\ref{Th:MTV_1}.
One generalization is given in an appendix of~\cite{MTV_Sh}, where it is conjectured that
for real parameters $\bs$, orbits of critical points of generalized master functions are real.
For the Lie algebra $\sln$, this holds as the polynomials $p_i$ of Section~\ref{S:construction}
are real. 
This new conjecture also holds for the Lie algebras $\mathfrak{sp}_{2n}$ and
$\mathfrak{so}_{2n+1}$, by the results in Section~7 of~\cite{MV04}. 

We also discuss other proofs of the Shapiro conjecture.

%
%
\subsection{Discrete Wronskians to Calogero-Moser spaces}\label{S:Calogero-Moser}

The XXX model is another integrable system studied by Mukhin, Tarasov, and Varchenko~\cite{MTV_XXX}.
This work is similar to their work on the periodic Gaudin model, including
Wronskians, a Bethe ansatz, and  symmetric operators.
One difference is that  $U\gln$ is replaced by the Yangians,
$Y\gln$, which are a deformation of the universal enveloping algebra of the current
algebra $\gln[t]$.
(The current algebra $\gln[t]$ consists of polynomials in $t$ with coefficients in $\gln$.) 
Another is that the usual Wronskian is replaced by the \DeCo{{\sl discrete Wronskian}},
which depends upon a real number $h$, 
\[
  \DeCo{\Wr_h(f_0,\dotsc,f_n)}\ :=\ \det\;\left(\begin{matrix}
   f_0(t)&f_0(t+h)&\dotsb&f_0(t+nh)\\
    f_1(t)&f_1(t+h)&\dotsb&f_n(t+nh)\\
    \vdots&\vdots&\ddots&\vdots\\
    f_n(t)&f_n(t+h)&\dotsb&f_n(t+nh)\\
  \end{matrix}\right)\ ,
\]
and the functions are \DeCo{{\sl quasi-polynomials}}, $f_i(t)=e^{b_i t}g_i(t)$, where  
$g_i(t)$ is a polynomial.
The linear span $V$ of quasi-polynomials
$e^{b_0 t}g_0(t), \dotsc, e^{b_n t}g_n(t)$ is a space of
quasi-polynomials.
This discrete Wronskian has the form
$w(t)e^{\sum b_i \cdot t}$, where $w(t)$ is a polynomial that is well-defined up
to a scalar.

\begin{theorem}[Theorem~2.1 of~\cite{MTV_XXX}]\label{Th:XXX}
  Let $V$ be a space of quasi-polynomials with discrete Wronskian
  $\Wr_h(V)=\prod_{i=1}^N(t-s_i)e^{\sum b_i\cdot t}$ whose roots are real and satisfy 
\[
   |s_i-s_j|\ \geq\ |h|\qquad\mbox{for all}\ i\neq j\,,
\]
 then $V$ has a basis of real quasi-polynomials.
\end{theorem}

This condition on the separation of roots cannot be relaxed if the theorem is to hold for all
exponents $b_i$.
When $n=1$ and $b_0=b_1=0$, this is a special case of the main theorem in Eremenko,
{\it et al.}~\cite{EGSV}.
We will not discuss the proof of Theorem~\ref{Th:XXX}, except to remark that it depends upon the
results of~\cite{MTV_06,MTV_07}.

In the limit as $h\to 0$, the discrete Wronskian becomes the usual Wronskian, which yields the
following theorem.

\begin{theorem}[Theorem~4.1 of~\cite{MTV_XXX}]\label{Th:XXX_cont}
  Let $V$ be a space of quasi-polynomials whose Wronskian has only real roots.
  Then $V$ has a basis of real quasi-polynomials.
\end{theorem}

When the exponents $b_i$ are all zero, this reduces to Theorem~\ref{Th:MTV_1}, and therefore
is a generalization of the Shapiro conjecture.
It implies Theorem~\ref{Th:matrix} from the Introduction.
Suppose that $b_0,\dotsc,b_n$ are distinct real numbers, $\alpha_0,\dotsc,\alpha_n$ are
complex numbers, and consider the matrix
 \begin{equation}\label{Eq:Z-form}
    Z\ :=\ \left(\begin{matrix}
      \alpha_0& (b_0-b_1)^{-1}& \dotsb & (b_0-b_n)^{-1}\\
      (b_1-b_0)^{-1}& \alpha_1& \dotsb & (b_1-b_n)^{-1}\\
         \vdots      & \vdots & \ddots & \vdots\\
      (b_n-b_0)^{-1}& (b_n-b_1)^{-1}& \dotsb & \alpha_n
    \end{matrix}\right)\ .
 \end{equation}
%

\noindent{\bf Theorem~\ref{Th:matrix}} (Theorem 5.4 of~\cite{MTV_XXX}){\bf .}
{\it 
  If $Z$ has only real eigenvalues, then the numbers $\alpha_0,\dotsc,\alpha_n$ are real.
}
\medskip

\noindent{\it Proof.}
 We follow~\cite{MTV_XXX}, deducing this from Theorem~\ref{Th:XXX_cont} and some matrix identities.
 Since
 \begin{equation}\label{Eq:derivative}
   \frac{d^m}{dt^m}(t-a)e^{bt}\ =\ b^m(t-a)e^{bt}\ +\ mb^{m-1}e^{bt}\,,
 \end{equation}
 if $A$ is the diagonal matrix $\diag(a_0,\dotsc,a_n)$, 
 $E:=\diag(e^{b_0t},\dotsc,e^{b_nt})$, $V$ is the Vandermonde matrix $(b_j^i)_{i,j=0}^n$, and 
 $W:=(ib_j^{i-1})_{i,j=0}^n$, then~\eqref{Eq:derivative} implies that
\[
   \left(  \frac{d^i}{dt^i}(t-a_j)e^{b_jt}\right)_{i,j=0}^n  \ =\ 
  \left[ V(It-A)\ +\ W\right] E\,,
\]
 and therefore
 \begin{multline}\label{Eq:Wr-formula}
   \qquad \Wr\bigl( (t-a_0)e^{b_0t}\,,\,\dotsc\,,\,(t-a_n)e^{b_nt}\bigr)\\
   \ =\ 
   e^{\sum b_i\cdot t}\cdot\prod_{i<j}(b_j-b_i)\cdot\det \left[It -(A - V^{-1}W)\right]\ .\qquad 
 \end{multline}

 We deduce a formula for $V^{-1}W$.
 The inverse of the Vandermonde comes from Lagrange's interpolation formula.
 For each $i=0,\dotsc,n$, set
\[
   \DeCo{\ell_i(u)}\ :=\ 
    \frac{\prod_{k\neq i}(u-b_k)}{\prod_{k\neq i}(b_i{-}b_k)}\ =\ 
   \sum_{j=0}^n \DeCo{\ell_{i,j}} u^j\ .
\]
 Since $\ell_i(b_j)=\delta_{i,j}$, we see that $V^{-1}=(\ell_{i,j})_{i,j=0}^n$.
 But then
\[
  V^{-1}W\ =\ (\ell'_i(b_j))_{i,j=0}^n\,.
\]
 This gives formulas for the entries of $V^{-1}W=(m_{i,j})_{i,j=0}^n$,
 \begin{eqnarray*}
  m_{i,j} &=& \frac{\prod_{k\neq i,j}(b_j-b_k)}{\prod_{k\neq i}(b_i-b_k)}\qquad
    \mbox{if}\quad i\neq j\quad\mbox{and}\\
  m_{i,i} &=& \sum_{k\neq i} \frac{1}{b_i-b_k}\ .
 \end{eqnarray*}

 Let $B$ be the diagonal matrix
 $\diag(\prod_{k\neq i}(b_i-b_k), i=0,\dotsc,n)$ and 
 $M$ be the diagonal of $V^{-1}W$.
 We leave the following calculation to the reader,
\[
   B^{-1}Z B\ =\ \diag(\alpha_0,\dotsc,\alpha_n)\ +\ M\ -\ V^{-1}W\,.
\]
 Combining this with~\eqref{Eq:Wr-formula}, we see that if 
 \begin{equation}\label{Eq:relation}
    a_i\ =\ \alpha_i\ +\ m_{i,i}\qquad i=0,\dotsc,n\,,
 \end{equation}
 then the eigenvalues of $Z$ are exactly the roots of the Wronskian~\eqref{Eq:Wr-formula}.

 Since the matrix $Z$ has only real eigenvalues, Theorem~\ref{Th:XXX_cont} implies that the
 span of $(t-a_0)e^{b_0t},\dotsc,(t-a_n)e^{b_nt}$ has a basis of real
 quasi-polynomials.
 Since the exponents $b_i$ are real and distinct, the numbers $a_0,\dotsc,a_n$ are real as are
 the entries of $V^{-1}W$, and
 so~\eqref{Eq:relation} implies that the entries $\alpha_i$ of $Z$ are real.
\hfill\raisebox{-4.5pt}{\includegraphics[height=14pt]{pictures/qed.eps}}\smallskip

 Theorem~\ref{Th:matrix} has an interesting consequence.

\begin{corollary}\label{Cor:real_eigenvalues}
 Suppose that $X$ and $Z$ are square complex matrices such that 
 \begin{equation}\label{Eq:CM_commutator}
   [X,Z]\ =\ I\ -\ K\,,
 \end{equation}
 where $K$ has rank $1$.
 If both $X$ and $Z$ have real eigenvalues, then they may be simultaneously conjugated to real
 matrices. 
\end{corollary}

\noindent{\it Proof.}
 It suffices to show this for a dense open subset of such pairs $(X,Z)$ of matrices.
 Suppose that $X$ is diagonalizable with eigenvalues $b_0,\dotsc,b_n$ and that we have
 conjugated $(X,Z)$ so that $X$ is diagonal.
 If we write $Z=(z_{i,j})_{i,j=0}^n$, then 
 \begin{equation}\label{EQ:rkOne}
   [X,Z]\ =\ \left( z_{i,j}(b_j-b_i)\right)_{i,j=0}^n\,.
 \end{equation}
 The rank 1 matrix $K$ has the form $(\beta_i\gamma_j)_{i,j=0}^n$, where $\beta,\gamma$ are
 complex vectors.
 By~\eqref{Eq:CM_commutator} and~\eqref{EQ:rkOne}, the diagonal entries of $K$ are all 1,
 so that $\beta_i\gamma_i=1$, so in fact $\beta,\gamma\in (\C^\times)^{n+1}$ with 
 $\gamma=\beta^{-1}$. 
 Conjugating~\eqref{Eq:CM_commutator} by $\beta$ (considered as a diagonal matrix), 
 we may assume that $K$ is the matrix whose every entry is 1, and so 
\[
   z_{i,j}(b_j-b_i)\ =\ \delta_{i,j}\ -\ 1\,,
\]
 or, if $i\neq j$, $b_i\neq b_j$, and $z_{i,j}=(b_i-b_j)^{-1}$.
 But then $Z$ has the form~\eqref{Eq:Z-form} (where $\alpha_i=z_{i,i}$), and
 Theorem~\ref{Th:matrix} implies that all of its entries are real.
\hfill\raisebox{-4.5pt}{\includegraphics[height=14pt]{pictures/qed.eps}}\vspace{5pt}

Mukhin, Tarasov, and Varchenko noted that Theorem~\ref{Th:XXX_cont} follows from
Theorem~\ref{Th:matrix} by the duality studied in~\cite{MTV_06b}, and that the Shapiro
conjecture for Grassmannians is the case of Corollary~\ref{Cor:real_eigenvalues} when $Z$ is
nilpotent. 
We close this section with an interesting circle of ideas related to
Corollary~\ref{Cor:real_eigenvalues}. 

Let $\overline{C}_n$ be the set of all pairs $(X,Z)$ of $(n{+}1)\times(n{+}1)$ complex matrices 
such that  $[X,Z]-I$ has rank 1.
The group $Gl_{n+1}(\C)$ acts on $\overline{C}_n$ by simultaneous conjugation and
Wilson~\cite{Wi98} defines the \DeCo{{\sl Calogero-Moser space $C_n$}} to be the quotient of
$\overline{C}_n$ by $Gl_{n+1}(\C)$.
He shows this is a smooth affine variety of dimension $2n$.
It has many incarnations.
It is the phase space of the (complex) Calogero-Moser integrable system~\cite{KKS},
Etingof and Ginzburg~\cite{EtGi02} showed that $C_n$ parametrizes irreducible representations
of a 
certain rational Cherednik algebra, and  Wilson's adelic Grassmannian~\cite{Wi98} is naturally
the union of all the spaces $C_n$.

Let $C_n(\R)$ be the real points of $C_n$.
This turns out to be image of the real points of $\overline{C}_n$ under the quotient map 
$\DeCo{\pi_n}\colon\overline{C}_n\to C_n$.
The map that takes a matrix to its eigenvalues induces a map
\[
    \DeCo{\Upsilon}\ \colon\  C_n\ \longrightarrow\ 
     \C^{(n+1)}\times \C^{(n+1)}\,,
\]
where $\C^{(n+1)}:=\C^{n+1}/\calS_{n+1}$.
Etingof and Ginzburg showed that $\Upsilon$ is a finite map of degree $(n{+}1)!~$
We restate Corollary~\ref{Cor:real_eigenvalues}.
\medskip

\noindent{\bf Corollary~\ref{Cor:real_eigenvalues}.} \quad
{\it $\Upsilon^{-1}(\R^{(n+1)}\times \R^{(n+1)})\subset C_n(\R)$.}\medskip

This in turn implies the Shapiro conjecture for Grassmannians, which is the case of
Corollary~\ref{Cor:real_eigenvalues} when $Z$ is nilpotent, 
\[
   \Upsilon^{-1}( \R^{(n+1)}\times\{0\})\ \subset\ C_n(\R)\,.
\]

The \DeCo{{\sl rational Cherednik algebra $H_n$}}~\cite{EtGi02} is generated by 
the polynomial subalgebras $\C[x_0,\dotsc,x_n]$ and $\C[z_0,\dotsc,z_n]$ and the group algebra
$\C\calS_{n+1}$ subject to 
 \begin{equation}\label{Eq:RCA_rels}
   \begin{array}{rclcrcl}
   \sigma_{ij} x_i&=&x_j\sigma_{ij} &\qquad& [x_i,z_j]&=& \sigma_{ij}
          \qquad\mbox{if}\quad i\neq j\,,\\ \rule{0pt}{14pt}
   \sigma_{ij} z_i&=&z_j\sigma_{ij} &\qquad& [x_i,z_i]&=& -\sum_{j\neq i}\sigma_{ij}\,,
 \end{array}
\end{equation}
where $\sigma_{ij}\in\calS_{n+1}$ is the transposition $(i,j)$.
The symmetrizing idempotent is
\[
    \DeCo{e}\ :=\ \frac{1}{(n{+}1)!}\sum_{\sigma\in\calS_{n+1}} \sigma\ .
\]
For $p\in C_n$, write \DeCo{$\C_p$} for the 1-dimensional representation of the coordinate ring 
of $C_n$ in which a function $f$ acts by the scalar $f(p)$.

\begin{theorem}[\cite{EtGi02}, Theorems 1.23 and 1.24]
 \ 
 \begin{enumerate}
   \item $eH_ne$ is isomorphic to the coordinate ring of $C_n$.

    \item Irreducible representations of $H_n$ are parametrized by the points $p$ of
      $C_n$, where the corresponding representation is
\[
    \DeCo{M_p}\ :=\ H_ne\otimes_{eH_ne} \C_p\,.
\]
  
    \item $M_p$ is isomorphic to $\C\calS_{n+1}$ as an $\calS_{n+1}$-module.
 \end{enumerate}
\end{theorem}

Etingof and Ginzburg connect the structure of the representations $M_p$ to the Calogero-Moser 
space.
Let $\calS_n$ act on the indices $\{0,\dotsc,n{-}1\}$. 
Then $x_n$ and $z_n$ both stabilize the subspace $M_p^{\calS_n}$ of invariants, which has
dimension $n{+}1$.

\begin{theorem}[\cite{EtGi02}, Theorem 11.16]\label{Th:CM_action}
  In any basis of $M_p^{\calS_n}$, $x_n,z_n$ act by a pair of matrices $(X,Z)\in\overline{C}_n$
  such that $\pi_n(X,Z)=p$.
\end{theorem}

%
%
\subsection{Transversality in the Shapiro conjecture}

While Mukhin, Tarasov, and Varchenko prove Theorem~\ref{Th:strong} in~\cite{MTV_R}, their actual
result is much deeper.
Each ramification sequence $\ba$ for $\Gr$ corresponds to a dominant weight
$\mu(\ba)$ for $\sln$ such that, given a Schubert problem $\bA:=(\ba^{(1)},\dotsc,\ba^{(m)})$, 
the intersection number $\delta(\bA)$ is equal to the dimension of the space of singular
vectors
\[
   \left(V_{\mu(\sba^{(1)})}\otimes V_{\mu(\sba^{(2)})}\otimes
    \dotsb \otimes V_{\mu(\sba^{(m)})}\right)[0]\,,
\]
as both numbers are determined by the same formula based on the 
Littlewood-Richardson rule.
The result of~\cite{MTV_R} links more subtle scheme-theoretic information about the
intersection of Schubert varieties to algebraic information about the action of commuting
operators on the singular vectors.

The coordinate ring $R_\bA(\bs)$ of an intersection of Schubert varieties~\eqref{Eq:tr_int}
is an Artin algebra of dimension $\delta(\bA)$, because the Pl\"ucker formula~\eqref{Eq:Pl} 
forces the intersection to be zero dimensional.
It is semisimple exactly when the intersection is transverse. 
Because of the Pl\"ucker formula, the intersection lies in the big Schubert cell 
$\Gr^\circ$, and so $R_\bA(\bs)$ is a quotient of the coordinate ring $\calR$ of $\Gr^\circ$.
Then the coregular representation of $R_\bA(\bs)$ on its dual $R_\bA(\bs)^*$ induces an action
of $\calR$ on $R_\bA(\bs)^*$.
This is the scheme-theoretic information.

Given a finite-dimensional representation $V$ of $\gln$ and a complex number $s$, requiring
$t$ to act on $V$ via scalar multiplication by $s$, gives the 
\DeCo{{\sl evaluation module $V(s)$}} of the current algebra $\gln[t]$.
If $\bmu=(\mu^1,\dotsc,\mu^m)$ are dominant weights and $\bs=(s_1,\dotsc,s_m)$ are distinct
complex numbers, then the evaluation module
\[
   \DeCo{V_\sbmu(\bs)}\ :=\ V_{\mu^1}(s_1)\otimes V_{\mu^2}(s_2)\otimes
           \dotsb\otimes V_{\mu^m}(s_m)
\]
is an irreducible $\gln[t]$-module~\cite{CP}.

The universal enveloping algebra $U\gln[t]$ has a commutative subalgebra \DeCo{$\calB$}, called
the \DeCo{{\sl Bethe algebra}}~\cite{Ga,Ta}. 
As $\calB$ commutes with $\gln$, it acts on spaces of singular vectors
in the evaluation module $V_{\sbmu}(\bs)$.
The action of the Bethe algebra on the singular vectors $V_{\sbmu}(\bs)[0]$ is the algebraic
information from representation theory. 
Let $\calB_\sbmu(\bs)$ be the image of $\calB$ in the endomorphism algebra
of $V_{\sbmu}(\bs)[0]$.

A main result in~\cite{MTV_R} is that these two actions, 
$\calR$ on $R_\bA(\bs)$ and $\calB$ on $V_{\sbmu}(\bs)[0]$, 
are isomorphic when $\bmu=(\mu(\ba^1),\dotsc,\mu(\ba^m))$, which we write as $\bmu(\bA)$.
This requires that we identify $\calR$ with $\calB$ in some way.
For that, the big cell $\Gr^\circ$ can be identified with $(n{+}1)\times (d{-}n)$ matrices
$(y_{i,j})$, whose entries identify $\calR$ with $\C[y_{i,j}]$.
The Bethe algebra has generators $B_{i,j}$, where $1\leq i\leq n{+}1$ and 
$1\leq j$.
Define the map $\tau\colon\calR\hookrightarrow\calB$ by 
\[
   \tau(y_{i,j})\ =\ B_{i,j}\,.
\]
Mukhin, Tarasov, and Varchenko also give a linear bijection
$\phi\colon R_\bA(\bs)^*\to V_{\sbmu}(\bs)[0]$.

\begin{theorem}
 Let $\bA=(\ba^1,\dotsc,\ba^m)$ be a Schubert problem for $\Gr$ and $\bs=(s_1,\dotsc,s_m)$ be
 distinct complex numbers.
 Then the map $\tau$ induces an isomorphism of algebras
 $\tau\colon R_\bA(\bs)\xrightarrow{\,\sim\,}\calB_{\sbmu(\bA)}(\bs)$ so that, for 
 $f\in R_\bA(\bs)$ and $g\in R_\bA(\bs)^*$,
 $\mu(f^* g)=\tau(f)\mu(g)$.
 That is, $(\tau,\mu)$ gives an isomorphism between the 
 coregular representation of $R_\bA(\bs)$ on its linear dual  $R_\bA(\bs)^*$ and the
 action of the Bethe algebra $B_{\sbmu(\bA)}(\bs)$ on the singular vectors $V_{\sbmu}(\bs)[0]$.
\end{theorem}

Theorem~\ref{Th:strong} now follows, as the image of the Bethe algebra on the
singular vectors $V_{\sbmu}(\bs)[0]$ is generated by the Gaudin Hamiltonians, which are
diagonalizable when the parameters $s_i$ are real.
Thus $\calB_{\sbmu(\bA)}(\bs)$ and hence $R_\bA(\bs)$ are semisimple, which implies
that the intersection of Schubert varieties~\eqref{Eq:tr_int} was transverse.

We remark that this uses the general form of the results in~\cite{MTV_Sh} which we did not
describe in these notes.
Also, the coincidence of numbers, $\delta(\bA)=\dim(V_{\sbmu(\bA)}(\bs)[0])$, is an important
ingredient in the proof that $\mu$ is a bijection.

Very recently, Mukhin, Tarasov, and Varchenko related this Bethe algebra to the center of the
rational Cherednik algebra of Section~\ref{S:Calogero-Moser}, and to the algebra of regular
functions on the Calogero-Moser space~\cite{MTV_CM}.

%
%
\section{Applications of the Shapiro conjecture}\label{S:applications}

Theorem~\ref{Th:MTV_1} and its stronger version, Theorem~\ref{Th:strong}, have 
other applications in mathematics.
Some are straightforward, such as linear series on $\PP^1$ with real ramification.
Here, we discuss this application in the form of maximally inflected curves and rational
functions with real critical points, as well as Purbhoo's considerably deeper application
in which he recovers the basic combinatorial algorithms of Young tableaux from the
monodromy of the Wronski map.
We close with Eremenko and Gabrielov's  computation of the degree of the real Wronski
map, which implies that when $d$ is even and $W$ is a generic real polynomial, there are many
spaces of real polynomials with Wronskian $W$.

%
%
\subsection{Maximally inflected curves}
One of the earliest occurrences of the central mathematical object of these notes, spaces
of polynomials with prescribed ramification, was in algebraic geometry, as these spaces of
polynomials are linear series $P\subset H^0(\PP^1,{\mathcal O}(d))$ on $\PP^1$ with
prescribed ramification.
Their connection to Schubert calculus originated in work of Castelnuovo in
1889~\cite{Ca1889} on $g$-nodal rational curves, and this was important in Brill-Noether
theory (see Ch.~5 of~\cite{HaMo}) and the Eisenbud-Harris theory of limit linear
series~\cite{EH83,EH87}. 

A linear series $P$ on $\PP^1$ of degree $d$ and dimension $n{+}1$ (a point in $\Gr$)
gives rise to a degree $d$ map
 \begin{equation}\label{Eq:rational_curve}
   \varphi\ \colon\ \PP^1\ \longrightarrow\ \PP^n\ =\ \PP(P^*)
 \end{equation}
of $\PP^1$ to projective space.
We will call this map a curve.
The linear series is \DeCo{{\sl ramified}} at points $s\in\PP^1$ where $\varphi$ is not
convex (the derivatives $\varphi(s),\varphi'(s),\dotsc,\varphi^{(n)}(s)$ do not span $\PP^n$). 
Call such a point $s$ a \DeCo{{\sl flex}} of the curve~\eqref{Eq:rational_curve}.

A curve is real when $P$ is real.
It is \DeCo{{\sl maximally inflected}} if it is real and all of its flexes are real.
The study of these curves was initiated in~\cite{KS}, where restrictions on the topology
of plane maximally inflected curves were established.

Let us look at some examples.
There are two types of real rational cubic curves, which are distinguished by their
singular points.
The singular point of the curve on the left below is a node and it is connected to the rest of
the curve, while the singular point on the other curve is isolated from the
rest of the curve, and is called a \DeCo{{\sl solitary point}}. 
\[
  \begin{picture}(67,85)
   \put(5,15){\includegraphics[height=70pt]{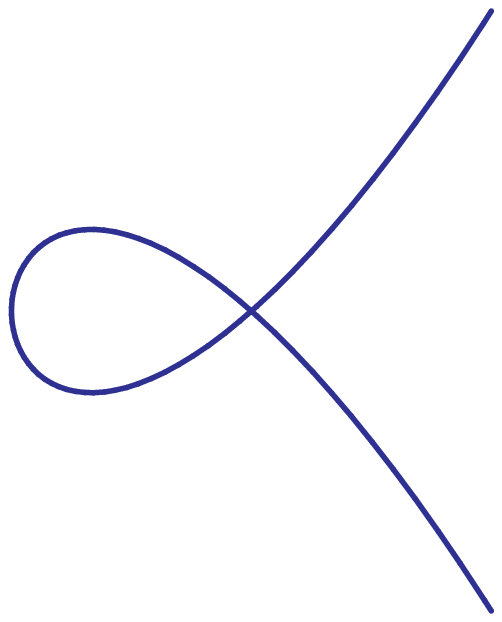}}
   \put(4,0){$y^2=x^3+x^2$}
  \end{picture}
  \qquad\qquad\qquad
  \begin{picture}(70,85)
   \put(0,15){\includegraphics[height=70pt]{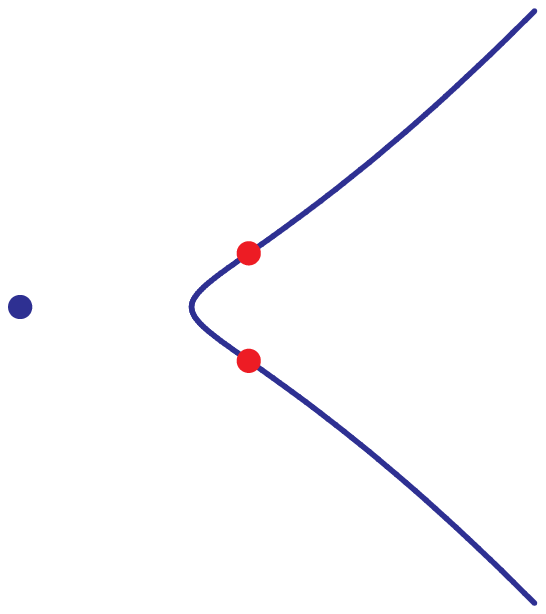}}
   \put(8,0){$y^2=x^3-x^2$}
  \end{picture}
\]
While both curves have one of their three flexes at infinity, only the curve on
the right has its other two flexes real (the dots) and is therefore maximally inflected.
A nodal cubic cannot be maximally inflected.

Similarly, a maximally inflected quartic with six flexes has either 1 or 0 of its
(necessarily 3) singular points a node, and necessarily 2 or 3 solitary points.
We draw the two types of maximally inflected quartics having six flexes, omitting their
solitary points.
Due to the symmetry in the placement of the flexes, the first quartic has two flexes on
its node---one for each branch through the node.
\[
  \includegraphics[height=80pt]{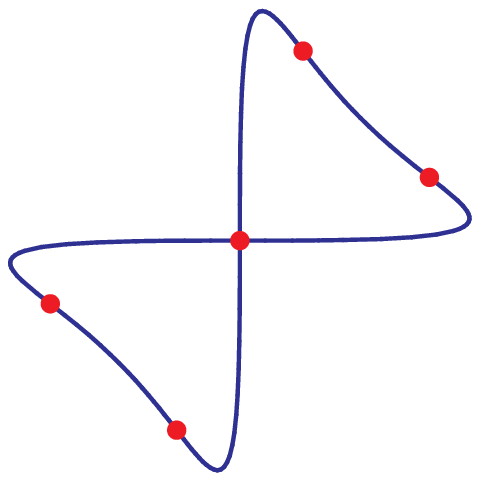}
   \qquad\qquad
  \includegraphics[height=80pt]{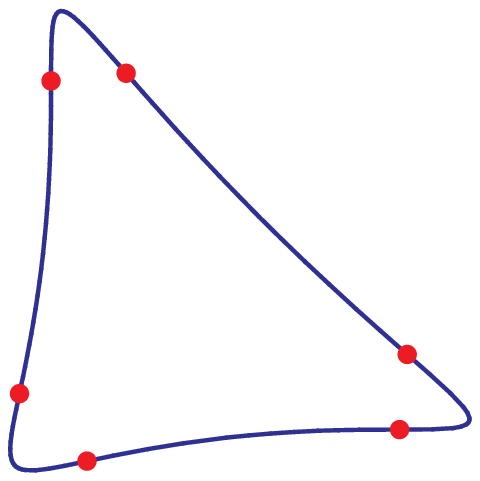}
\]

When $d\geq 3$, the image of a rational curve $\gamma\colon\PP^1\to\PP^2$ of degree $d$ is a
singular curve of arithmetic genus $\binom{d-1}{2}$.
In general, the singularities consist of $\binom{d-1}{2}$  ordinary double points, which
are where two smooth branches of the curve meet transversally.
A real rational curve has three types of such double points.
We have already seen nodes and solitary points.
The third kind of real double point consists of a pair of complex conjugate double points,
and is invisible in $\R\PP^2$.

The examples we gave had few nodes.
This is always the case.

\begin{theorem}[Corollary~3.3 and Theorem~4.1 of~\cite{KS}]\label{Th:topology}
  Given a maximally inflected plane curve of degree $d$, let $\delta,\eta,c$ be its
  numbers of solitary points, nodes, and pairs of complex conjugate double points.
  Then we have
\[
    d-2\ \leq\ \delta\ \leq \binom{n-1}{2}
    \qquad\mbox{and}\qquad 
   0\ \leq\ \eta+2c\ \leq\ \binom{n-2}{2}\,.
\]
  Furthermore, there exist maximally inflected curves of degree $d$ 
  with $\binom{n-1}{2}$ solitary points (and hence no other singularities), and 
  there exist curves with $\binom{n-2}{2}$ nodes and $d{-}2$ solitary points.
\end{theorem}

While many constructions of maximally inflected curves were known, Theorem~\ref{Th:MTV_1}, 
and in particular Theorem~\ref{Th:strong}, show that there are many maximally inflected
curves:
Any curve $\varphi\colon\PP^1\to\PP^n$ whose ramification lies in $\R\PP^1$ must be real
and is therefore maximally inflected.

Theorem~\ref{Th:topology} is proven using the Pl\"ucker formula~\eqref{Eq:Pl} and the Klein
formula from topology.
This result and some intriguing calculations in Section 6 of~\cite{KS} suggest that the
number of solitary points is a deformation invariant.
That is, if the points of inflection move, then the number of solitary points is constant.
Examples show that the number of nodes may change under a continuous deformation of the
inflection points, with a pair of nodes colliding to become a pair of complex conjugate double
points, but we have not observed collisions of solitary points.

%
%
\subsection{Rational functions with real critical points}
A special case of Theorem~\ref{Th:MTV_1}, proved earlier by
Eremenko and Gabrielov, serves to illustrate the breadth of mathematical areas touched by
this Shapiro conjecture.
When $n=1$, we may associate a rational function $\varphi_P:=f_1(t)/f_2(t)$ 
to a basis $\{f_1(t),f_2(t)\}$ of a vector space $P\in\G(1,d)$ of polynomials.
Different bases give different rational functions, but they all differ from each other by
a fractional linear transformation of the image $\PP^1$.
We say that such rational functions are \DeCo{{\sl equivalent}}.

The critical points of a rational function are the points of the
domain $\PP^1$ where the derivative of $\varphi_P$, 
\[
   d \varphi_P\ :=\ 
   \frac{f'_1 f_2 - f_1 f_2'}{f_2^2}\ =\ 
   \frac{1}{f_2^2}\cdot \det
     \left(\begin{matrix}f_1&f_2\\f'_1&f'_2\end{matrix}\right)\ ,
\]
vanishes.
That is, at the roots of the Wronskian.
These critical points only depend upon the equivalence class of the rational function.
Eremenko and Gabrielov~\cite{EG02a} prove the following result about the critical points
of rational functions.

\begin{theorem}\label{Th:EG_rat}
  A rational function $\varphi$ whose critical points lie on a circle in $\PP^1$ maps that
  circle to a circle.
\end{theorem}

To see that this is equivalent to Theorem~\ref{Th:MTV_1} when $n=1$, note that we may apply a
change of variables to $\varphi$ so that its critical points lie on the circle
$\R\PP^1\subset\PP^1$. 
Similarly, the image circle may be assumed to be $\R\PP^1$.
Reversing these coordinate changes establishes the equivalence.

The proof used methods specific to rational functions.
Goldberg showed~\cite{Go91} that there are at most $\DeCo{c_d}:=\frac{1}{d}\binom{2d-2}{d-1}$
rational functions of degree $d$ with a given collection of $2d-2$ simple critical points.
If the critical points of a rational function $\varphi$ of degree $d$ lie on a circle
$C\subset\C\PP^1$ and if $\varphi$ maps $C$ to $C$, then $\varphi^{-1}(C)$ forms a graph on
the Riemann sphere with nodes the $2d{-}2$ critical points, each of degree $4$, and each
having two edges along $C$ and one edge on each side of $C$.
We mark one of the critical points to fix the ordered list of the critical
points.
It turns out that there are also $c_d$ such abstract graphs with a distinguished vertex.
Call these graphs \DeCo{{\sl nets}}.
(In fact, $c_d$ is Catalan number, which counts many objects in 
combinatorics~\cite[Exer.~6.19, p.~219]{Sta}.)
Here are the $c_4=5$ nets for $d=4$:
\[
   \includegraphics{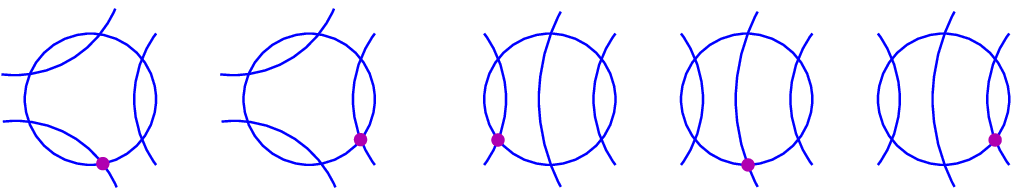}
\]
Eremenko and Gabrielov used the uniformization theorem in complex analysis to construct such a
rational function $\varphi$ for each net and choice of critical points on $C$. 
Since $c_d$ is the upper bound for the number of such rational functions, this
gave all rational functions with given set of critical points and thus
proved Theorem~\ref{Th:EG_rat}.

More recently, Eremenko and Gabrielov~\cite{EG05} found an elementary proof of
Theorem~\ref{Th:EG_rat} which is based upon analytic continuation and a very refined version of
the construction underlying Theorem~\ref{Th:local}.
This has unfortunately never been published.

By Theorem~\ref{Th:local}, there exists a Wronski polynomial $W_0(t)$ of degree $2d{-}2$ with
distinct real roots for which there are $c_d$ spaces of real polynomials with Wronskian
$W_0(t)$. 
Suppose that $W_0$ is a member of a continuous family $W_\tau$ for $\tau\in[0,1]$ of
polynomials of degree $2d{-}2$ with distinct real roots.
Since there are $c_d$ distinct spaces of polynomials with Wronskian $W_0$, there are $c_d$
distinct lifts of the path $W_\tau$ to paths of spaces of polynomials with Wronskian $W_\tau$,
at least for $\tau$ near zero.
The obstruction to analytically continuing these $c_d$ lifts occurs at critical points of the
Wronski map $\Wr\colon\G(1,d)_\R\to\R\PP^{2d-2}$.
Since this map is at most $c_d$ to 1, the first critical point in a fiber is a point where two 
of the lifted paths collide.

Eremenko and Gabrielov show that such a collision cannot occur.
The reason is simple: nets are constant along paths of spaces of polynomials in $\G(1,d)$ whose
Wronskian has $2d-2$ distinct roots, and each of the spaces of polynomials above $W_0(t)$ has
a different net.
Thus each lifted path has a different net, and no collision is possible.
They show that nets are constant along paths by a simple set-theoretic/topological argument.
A similar elementary argument applied to the construction of the spaces of polynomials in
Theorem~\ref{Th:local} shows that each space has a distinct net.
The proof is completed by observing that any Wronski polynomial may be joined to $W_0$
along some path $W_\tau$.

The elementary and constructive nature of this proof suggests that the Shapiro Conjecture for
Grassmannians may have an elementary proof, if a suitable substitute for nets can be
found when $n>1$.

%
%
\subsection{Tableaux combinatorics}
Starting from Theorems~\ref{Th:MTV_1} and~\ref{Th:strong} for the Schubert
problem $\bi_{n,d}$, Purbhoo~\cite{Purbhoo} shows that the basic combinatorial
properties and algorithms for Young tableaux are realized geometrically via the monodromy
groupoid of the Wronski map, $\Wr\colon\Gr\to\PP(\C_{(n+1)(d-n)}[t])$. 
In particular, Sch\"utzenberger slides, evacuation, Knuth equivalence and dual equivalence 
all arise geometrically.
Purbhoo uses his analysis of the monodromy groupoid to get a new proof of the
Littlewood-Richardson rule. 

A \DeCo{{\sl partition}} is a weakly decreasing sequence of nonnegative integers
$\blambda\colon\lambda_0\geq \lambda_1\geq\dotsb\geq\lambda_n\geq 0$.
We impose the restriction that $n{-}d\geq \lambda_0$.
Partitions are ramification sequences in disguise, with 
$\ba\colon 0\leq a_0<a_1<\dotsb<a_n\leq d$ corresponding to 
 \begin{equation}\label{Eq:lambda(a)}
  \DeCo{\blambda(\ba)}\ \colon\ a_n-n\ \geq\ \dotsb\ \geq\ a_1-1\ \geq\ a_0\,.
 \end{equation}

We identify a partition with its diagram, which is a left-justified array of boxes with
$\lambda_i$ boxes in the $i$th row.
For example,
\[
   \blambda=5322\ \longleftrightarrow\ \raisebox{-12pt}{\includegraphics{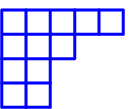}}\ .
\]
Write $|\blambda|$ for the number of boxes in $\blambda$.
By~\eqref{Eq:lambda(a)}, $|\ba|=|\blambda(\ba)|$.

The partial order on ramification sequences induces the partial order of componentwise
comparison on partitions, which is the inclusion of their diagrams.
The minimal partition is $\emptyset$ and the maximal partition (for us) is
$(d{-}n,\dotsc,d{-}n)$, which has $d{-}n$ repeated $n{+}1$ times.
Write this as $\maxp$.
Given $\bmu\leq\blambda$, the \DeCo{{\sl skew partition} $\blambda/\bmu$} is the difference of
their diagrams.
We set $\DeCo{|\blambda/\bmu|}:=|\blambda|-|\bmu|$.
For example,
\[
   5322/21\ =\ 
    \raisebox{-12pt}{\includegraphics{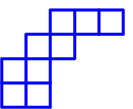}}
  \qquad\mbox{and}\qquad
  |5322/21|\ =\ 9\,.
\]

A \DeCo{{\sl standard Young tableau}} of \DeCo{{\sl shape}} $\blambda/\bmu$ is a filling of the
boxes of $\blambda/\bmu$ with the integers $1,2,\dotsc,|\blambda/\bmu|$ so that the entries increase
across each row and down each column.
Here are three fillings of the shape $331/1$; only the first two are tableaux.
{\small\[
  \begin{picture}(37,37)(-3.5,-3.5)
   \put(-3.5,-3.5){\includegraphics{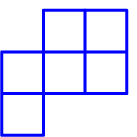}}
                    \put(12,24){$1$} \put(24,24){$3$}
   \put(0,12){$2$}  \put(12,12){$5$} \put(24,12){$6$} 
   \put(0, 0){$4$}  
  \end{picture}
    \qquad
  \begin{picture}(37,37)(-3.5,-3.5)
   \put(-3.5,-3.5){\includegraphics{pictures/331_1.eps}}
                    \put(12,24){$2$} \put(24,24){$4$}
   \put(0,12){$1$}  \put(12,12){$3$} \put(24,12){$5$} 
   \put(0, 0){$6$}  
  \end{picture}
    \qquad
  \begin{picture}(37,37)(-3.5,-3.5)
   \put(-3.5,-3.5){\includegraphics{pictures/331_1.eps}}
                    \put(12,24){$2$} \put(24,24){$6$}
   \put(0,12){$5$}  \put(12,12){$4$} \put(24,12){$3$} 
   \put(0, 0){$1$}  
  \end{picture}
\]}
Let $\SYT(\blambda/\bmu)$ be the set of all standard Young tableaux of shape $\blambda/\bmu$.

The degree $\delta(\bi_{n,d})$ of the Wronski map equals the cardinality of $\SYT(\maxp)$.
By Theorem~\ref{Th:strong}, the Wronski map is unramified over the locus of polynomials
with distinct real roots, and so the points in each fiber are in bijection with the set 
\DeCo{$\SYT(\maxp)$}.
This identification is almost canonical because in the region where the roots of the Wronskian
are clustered~\eqref{Eq:cluster} the identification is canonical, by the work of Eremenko and
Gabrielov~\cite{EG02c}, and the Wronski map is unramified over the locus of polynomials with
distinct roots.
Since nets are in natural bijection with tableaux, this identification for $n=1$ was done
by Eremenko and Gabrielov in~\cite{EG05}.

This identification can be extended to skew tableaux.
Given a partition $\blambda$, its dual is 
$\blambda^\vee\colon d{-}n{-}\lambda_n\geq\dotsb\geq d{-}n{-}\lambda_0$.
For partitions $\bmu\leq\blambda$, set
\[
  \DeCo{\G(\blambda/\bmu)}\ :=\ \Omega_{\ba(\sbmu)}F_\bullet(0)\;\bigcap\;
         \Omega_{\ba(\sbl^\vee)}F_\bullet(\infty)\,.
\]
The Wronskian of a space of polynomials $P\in\G(\blambda/\bmu)$ has degree at most
$|\blambda|$ and vanishes to order least $|\bmu|$ at zero.
Let \DeCo{$\PP(\blambda/\bmu)$} be the projective space of such polynomials.
This has dimension $|\blambda/\bmu|$, which is equal to the dimension of $\G(\blambda/\bmu)$. 
The restriction of the Wronski map to  $\G(\blambda/\bmu)$,
\[
  \Wr\ \colon\ \G(\blambda/\bmu) \ \longrightarrow\ \PP(\blambda/\bmu)\,,
\]
is finite, flat, and has degree equal to the cardinality of $\SYT(\blambda/\bmu)$.
Lastly, the Wronski map is unramified over the locus of polynomials in $\PP(\blambda/\bmu)$
with $|\blambda/\bmu|$ distinct nonzero real roots,
and there is an identification of the fibers with $\SYT(\blambda/\bmu)$.

Purbhoo gives an explicit identification of the fibers of the Wronski map by
extending the notion of standard tableaux.
Let $\bs=\{s_1,\dotsc,s_{|\sbl/\sbmu|}\}\subset \R\PP^1$ be a set of $|\blambda/\bmu|$
real numbers, possibly including $\infty$, that satisfy
\begin{enumerate}
 \item[(I)]  If $i\neq j$, then $|s_i|\neq|s_j|$.
 \item[(II)] $0\in\bs$ only if $\bmu=\emptyset$ and 
             $\infty\in\bs$ only if $\blambda=\maxp$.
\end{enumerate}
We identify such a subset $\bs$ with the polynomial 
$\DeCo{W_\bs}:=t^{|\sbmu|}\prod_{s\in\bs}(t-s)$ in $\PP(\blambda/\bmu)$ vanishing at $\bs$.
A \DeCo{{\sl standard Young tableau}} of \DeCo{{\sl shape}} $\blambda/\bmu$ with entries in
$\bs$ is a filling of the boxes of $\blambda/\bmu$ with elements of $\bs$ such that if we
replace each entry $s_i$ with its absolute value $|s_i|$, then the entries increase across
each row and down each column.

Let \DeCo{$\SYT(\blambda/\bmu;\bs)$} be the set of all standard Young tableaux of shape
$\blambda/\bmu$ with entries in $\bs$.
Replacing each entry $s_i$ in a tableau by $k$ if $s_i$ has the $k$th smallest absolute
value in $\bs$ defines the map $\ord\colon\SYT(\blambda/\bmu;\bs)\to\SYT(\blambda/\bmu)$.
For example,
{\small
\newcommand{\ms}{\hspace{5pt}}
\[
  \begin{picture}(55,55)(-2.5,-6.5)
   \put(-2.5,-6.5){\includegraphics{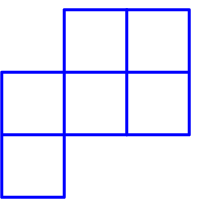}}
                      \put(18,36){$\sqrt{2}$} \put(36,36){\ms$4$}
   \put(0,18){\ms$e$} \put(18,18){$-8$}   \put(36,18){$\pi^2$} 
   \put(0, 0){$-6$}  
  \end{picture}
  \quad
    \raisebox{23pt}{\LARGE$\stackrel{\mbox{\small$\ord$}}{\longmapsto}$}
  \quad
  \begin{picture}(55,55)(-2.5,-6.5)
   \put(-2.5,-6.5){\includegraphics{pictures/331_1b.eps}}
                      \put(18,36){\ms$1$} \put(36,36){\ms$3$}
   \put(0,18){\ms$2$} \put(18,18){\ms$5$} \put(36,18){\ms$6$} 
   \put(0, 0){\ms$4$}  
  \end{picture}
\]}\vspace{-12pt}

Let $\bs(\tau)$ for $\tau\in[a,b]$ be a continuous path of subsets of
$\R\PP^1$ where $\bs(\tau)$  satisfies Conditions (I) and (II), except for finitely many points
$\tau\in (a,b)$ at which Condition (I) is violated exactly once in that
$s_i(\tau)=-s_j(\tau)$ for some $i\neq j$.
A path that is transverse to the locus where $s_i=-s_j$ for all $i\neq j$ is  \DeCo{{\sl generic}}.
Given a standard Young tableau $T_a$ of shape $\blambda/\bmu$ and filling $\bs(a)$, we can try
to lift $T_a$ to a family $T_\tau$ of standard tableaux for all $\tau\in[a,b]$.
We do this by requiring that the entries in the boxes of $T_\tau$ vary continuously, unless the  
condition that $T_\tau$ forms a tableau is violated.

In any interval where  $\bs(\tau)$ satisfies Condition (I), the entries of
$T_\tau$ vary continuously and $\ord(T_\tau)$ is constant.
Suppose that $\tau_0$ is a point of the path where Condition (I) is violated, and that
$s_i(\tau_0)=-s_j(\tau_0)$ is the pair witnessing this violation.
If $s_i$ and $s_j$ are in different rows and columns, they remain in their respective boxes
as $\tau$ passes $\tau_0$ and $\ord(T_\tau)$ changes as $\tau$ passes $\tau_0$.
If $s_i$ and $s_j$ are in the same row or column, then they are adjacent and 
leaving them in their respective boxes violates the condition that $T_\tau$ is a tableau,
so we require them to switch places and $\ord(T_\tau)$ does not change as $\tau$ passes
$\tau_0$. 

Given a generic path $\bs(\tau)$ for $\tau\in[a,b]$ and a tableau
$T_a\in\SYT(\blambda/\bmu;\bs(a))$, define \DeCo{$\slide_{\bs(\tau)}(T_a)$} to be the result
of this process applied to $T_a$.
This gives a bijection between $\SYT(\blambda/\bmu;\bs(a))$ and $\SYT(\blambda/\bmu;\bs(b))$.

\begin{example}
 We show this on a tableau of shape $(4,4,2)$, for the path
 $\bs(\tau)=\{\tau, -1, \dotsc,-9\}$ for $\tau\in[0,10]$.
 We only display when the tableau $T_\tau$ changes.
{\small
\newcommand{\tms}{\hspace{5pt}\DeCo{$\tau$}}
\[
  \begin{picture}(73,55)(-2.5,-6.5)
   \put(0,36){\includegraphics{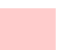}}
   \put(-2.5,-6.5){\includegraphics{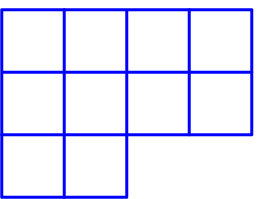}}
   \put(0,36){\tms} \put(18,36){$-1$} \put(36,36){$-3$} \put(54,36){$-8$}
   \put(0,18){$-2$} \put(18,18){$-4$} \put(36,18){$-6$} \put(54,18){$-9$} 
   \put(0, 0){$-5$} \put(18, 0){$-7$} 
  \end{picture}
  \ 
    \raisebox{23pt}{$\xrightarrow{\ \tau=1\ }$}
  \ 
  \begin{picture}(73,55)(-2.5,-6.5)
   \put(18,36){\includegraphics{pictures/cbox.eps}}
   \put(-2.5,-6.5){\includegraphics{pictures/442.eps}}
   \put(0,36){$-1$} \put(18,36){\tms} \put(36,36){$-3$} \put(54,36){$-8$}
   \put(0,18){$-2$} \put(18,18){$-4$} \put(36,18){$-6$} \put(54,18){$-9$}  
   \put(0, 0){$-5$} \put(18, 0){$-7$} 
  \end{picture}
  \ 
    \raisebox{23pt}{$\xrightarrow{\ \tau=3\ }$}
  \ 
  \begin{picture}(73,55)(-2.5,-6.5)
   \put(36,36){\includegraphics{pictures/cbox.eps}}
   \put(-2.5,-6.5){\includegraphics{pictures/442.eps}}
   \put(0,36){$-1$} \put(18,36){$-3$} \put(36,36){\tms} \put(54,36){$-8$}
   \put(0,18){$-2$} \put(18,18){$-4$} \put(36,18){$-6$} \put(54,18){$-9$}  
   \put(0, 0){$-5$} \put(18, 0){$-7$} 
  \end{picture}
  \ 
    \raisebox{23pt}{$\xrightarrow{\ \tau=6\ }$}
  \ 
\]
\[
  \begin{picture}(73,55)(-2.5,-6.5)
   \put(36,18){\includegraphics{pictures/cbox.eps}}
   \put(-2.5,-6.5){\includegraphics{pictures/442.eps}}
   \put(0,36){$-1$} \put(18,36){$-3$} \put(36,36){$-6$} \put(54,36){$-8$}
   \put(0,18){$-2$} \put(18,18){$-4$} \put(36,18){\tms} \put(54,18){$-9$}  
   \put(0, 0){$-5$} \put(18, 0){$-7$} 
  \end{picture}
  \ 
    \raisebox{23pt}{$\xrightarrow{\ \tau=9\ }$}
  \ 
  \begin{picture}(73,55)(-2.5,-6.5)
   \put(54,18){\includegraphics{pictures/cbox.eps}}
   \put(-2.5,-6.5){\includegraphics{pictures/442.eps}}
   \put(0,36){$-1$} \put(18,36){$-3$} \put(36,36){$-6$} \put(54,36){$-8$}
   \put(0,18){$-2$} \put(18,18){$-4$} \put(36,18){$-9$} \put(54,18){\tms}  
   \put(0, 0){$-5$} \put(18, 0){$-7$} 
  \end{picture}
  \ 
    \raisebox{23pt}{$\xrightarrow{\, \tau=10\,}$}
  \ 
  \begin{picture}(73,55)(-2.5,-6.5)
   \put(54,18){\includegraphics{pictures/cbox.eps}}
   \put(-2.5,-6.5){\includegraphics{pictures/442.eps}}
   \put(0,36){$-1$} \put(18,36){$-3$} \put(36,36){$-6$} \put(54,36){$-8$}
   \put(0,18){$-2$} \put(18,18){$-4$} \put(36,18){$-9$} \put(55.5,18){$10$}  
   \put(0, 0){$-5$} \put(18, 0){$-7$} 
  \end{picture}\hspace{33pt}
 \]}
 The combinatorial enthusiast will note that the box containing $\tau$ has just performed
 a Sch\"utzenberger slide through the subtableau formed by the negative entries.
 For comparison, we show the tableaux $\ord(T_0)$ and $\ord(T_{10})$.
{\small\[
  \begin{picture}(49,37)(-3.5,-3.5)
   \put(-3.5,-3.5){\includegraphics{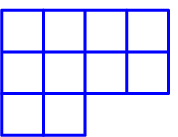}}
   \put(0,24){$1$} \put(12,24){$2$} \put(24,24){$4$} \put(36,24){$9$}
   \put(0,12){$3$} \put(12,12){$5$} \put(24,12){$7$} \put(34,12){$10$}
   \put(0, 0){$6$} \put(12, 0){$8$}   
  \end{picture}
  \qquad
  \begin{picture}(49,37)(-3.5,-3.5)
   \put(-3.5,-3.5){\includegraphics{pictures/442_s.eps}}
   \put(0,24){$1$} \put(12,24){$3$} \put(24,24){$6$} \put(36,24){$8$}
   \put(0,12){$2$} \put(12,12){$4$} \put(24,12){$9$} \put(34,12){$10$}
   \put(0, 0){$5$} \put(12, 0){$7$}   
  \end{picture}
\]}
\end{example}

We give Purbhoo's main theorem about the monodromy groupoid of
the Wronski map $\Wr\colon\G(\blambda/\bmu)\to\PP(\blambda/\bmu)$.

\begin{theorem}[\cite{Purbhoo}, Theorem 3.5]\label{Th:P_main}
 For each $\bs=\{s_1,\dotsc,s_{|\sbl/\sbmu|}\}\subset\R\PP^1$ satisfying Conditions {\rm (I)}
 and {\rm (II)}, there is a correspondence $P\leftrightarrow T(P)$ between points
 $P\in\G(\blambda/\bmu)$ with Wronskian $W_\bs$ and tableaux $T(P)\in\SYT(\blambda/\bmu;\bs)$.
 Under this correspondence, if $\bs(\tau)\subset\R\PP^1$ is a generic path for $\tau\in[a,b]$
 and $P_\tau$ is any lifting of that path to $\G(\blambda/\bmu)$, then 
\[
   T(P_b)\ =\ \slide_{\bs(\tau)} T(P_a)\ .
\]
\end{theorem}

Thus the combinatorial operation of sliding a tableau along a generic path $\bs(\tau)$
exactly describes analytic continuation in the fibers of the Wronski map above that path.
This sliding operation contains Sch\"{u}tzenberger's {\it jeu de taquin}~\cite{Schu1},
and much of tableaux combinatorics~\cite{Fu,Sagan,Sta} may be recovered from the geometry of
the Wronski map. 

Suppose $\bs=\{s_1,\dotsc,s_{|\sbl/\sbmu|}\}$ with $|s_1|<\dotsb<|s_{|\sbl/\sbmu|}|$.
If $T\in\SYT(\blambda/\bmu;\bs)$ and $\bt=\{s_i,\dotsc,s_j\}$ with $i<j$, then the entries
of $T$ in the set $\bt$ form a subtableau \DeCo{$T|_\bt$}.
Now suppose that $\bs=\bt\cup\bu$ where the elements of $\bt$ are positive, those of $\bu$
are negative, and we additionally have that $|t|<|u|$ for $t\in\bt$ and $u\in\bu$.
Write $|\bt|<|\bu|$ when this occurs.
Let $\bt'$ be a set of $|\blambda/\bmu|$ positive numbers with $|\bu|<|\bt'|$
and suppose that $\bs(\tau)$ for $\tau\in[0,1]$ is a generic 
path from $\bs=\bt\cup\bu$ to $\bs'=\bu\cup\bt'$.
Given a tableau $S\in\SYT(\blambda/\bmu;\bs)$, let
$S':=\slide_{\bs(\tau)}S\in\SYT(\blambda/\bmu;\bs')$, and define the subtableaux,
\[
   \DeCo{T}\ :=\ S|_\bt\,,\qquad
   \DeCo{U}\ :=\ S|_\bu\,,\qquad
   \DeCo{T'}\ :=\ S|_{\bt'}\,,\quad\mbox{and}\quad
   \DeCo{U'}\ :=\ S|_\bu\,.
\]
Because $|\bt|<|\bu|<|\bt'|$, $T$ is inside of $U$ and during the slide $T$ and $U$
move through each other to obtain the tableaux $U'$ and $T'$ with $U'$ inside of $T'$.
Schematically,
\[
  \begin{picture}(91,41)
   \put(0,0){\includegraphics{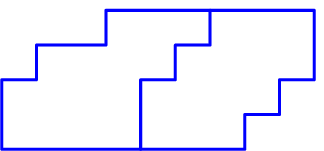}}
   \put(16,17){$T$}  \put(54,8){$U$}
  \end{picture}
   \qquad \raisebox{17pt}{$\xrightarrow{\ \mbox{$\slide_{\bs(\tau)}$}\ }$} \qquad
  \begin{picture}(91,41)
   \put(0,0){\includegraphics{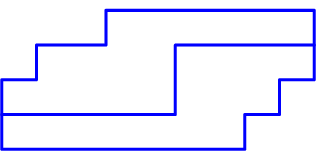}}
   \put(22,18){$U'$}  \put(60,14){$T'$}
  \end{picture}
\]
We write
\[
   U'\ =\ \slide_T U
   \qquad\mbox{and}\qquad
   T'\ =\ \slide_U T\,.
\]
Reversing the path $\bs(\tau)$ enables the definition of $U=\slide_{T'}U'$ and
 $T=\slide_{U'}T'$.
These notions are independent of the choice of path $\bs(\tau)$, by Theorem~\ref{Th:P_main}.
In fact $\slide_TU$ does not depend upon the set $\bt'$.
This geometrically defined operation was studied from a combinatorial perspective~\cite{BSS},
 where it was called tableaux switching, and its independence from choices was Theorem
 2.2(4)~{\em ibid}.

\begin{definition}
 Let $\bu$ be a set of negative numbers.
 Two tableaux $U_1,U_2\in\SYT(\blambda/\bmu;\bu)$ are \DeCo{{\sl equivalent}} if, for any set
 $\bt$ of $|\bmu|$ positive numbers with $|\bt|<|\bu|$ and any $T_1,T_2\in\SYT(\bmu;\bt)$, we
 have
\[
   \slide_{T_1}U_1\ =\ \slide_{T_2}U_2\,.
\]
 Two tableaux $U_1,U_2\in\SYT(\blambda/\bmu;\bu)$ are \DeCo{{\sl dual equivalent}} if, for any
 sets $\bt,\bt'$ of positive numbers with $|\bt|<|\bu|<|\bt'|$, shapes $\bmu/\bnu$,
 $\bk/\blambda$, and tableaux $T\in\SYT(\bmu/\bnu;\bt)$ and
 $T'\in\SYT(\bk/\blambda;\bt')$, each pair
 \[
    (\slide_T U_1\,,\; \slide_T U_2)\qquad\mbox{and}\qquad
    (\slide_{T'} U_1\,,\; \slide_{T'} U_2)
\]
 has the same shape.
 Replacing numbers by their negatives extends these definitions to tableaux with positive
 entries. 
\end{definition}

These are the fundamental equivalence relations on tableaux of Knuth-equivalence and of
Haiman's dual equivalence~\cite{Ha}. 
Purbhoo shows that these combinatorial equivalence relations coincide with
geometrically-defined relations that come from nonreduced fibers of the Wronski map. 

Suppose that $\bs=\{s_1,\dotsc,s_{(n+1)(d-n)}\}$ is a subset of $\R\PP^1$ satisfying
Condition (I) with $|s_1|<\dotsb<|s_{(n+1)(d-n)}|$, and suppose that 
$\bt=\{s_i,s_{i+1},\dotsc,s_j\}$ are the positive elements of $\bs$ and let $\bu=\bs-\bt$ be
its nonpositive elements.
Pick a positive number $a\in[s_i,s_j]$ and consider any path $\bs(\tau)$ for $\tau\in[0,1]$
that satisfies Condition (I) for $\tau\in[0,1)$ with $\bs(0)=\bs$, has constant nonpositive 
elements $\bu$, but whose positive elements all approach $a$ as $t\to 1$ so that
$\bs(1)=\{s_1,\dotsc,s_{i-1},a,\dotsc,a,s_{j+1},\dotsc,s_{(n+1)(d-n)}\}$.

Given a tableau $T\in\SYT(\maxp;\bs)$ corresponding to a point $P_T\in\Gr$ with Wronskian
$W_\bs$, we may analytically continue $P_T$ in the fibers of the Wronski map over the path
$W_{\bs(\tau)}$.
When $\tau\neq 1$, this continuation will be $P_{T_\tau}$, but when $\tau=1$ it will be
$\lim_{\tau\to 1}P_{T_\tau}$.
Write $P_T(\tau)$ for these points.
For each $\tau<1$ the points $P_T(\tau)$ are distinct for different
$T\in\SYT(\maxp;\bs(\tau))$, but in the limit as $\tau\to 1$ 
some paths may coalesce, as the fiber of the Wronskian is nonreduced at $\bs(1)$.

\begin{theorem}[\cite{Purbhoo}]
 Let $T,T'\in\SYT(\maxp;\bs)$.
 Then $T_\bt$ is equivalent to $T'|_{\bt}$ if and only if 
 $P_T(1)=P_{T'}(1)$.
\end{theorem}

Let $\bs'(\tau)$ for $\tau\in[0,1]$ be another generic path with
$\bs'(0)=\bs$ in which the positive elements are constant, but the others 
converge to some fixed negative number $a$.
We define $P'_T(\tau)$ to be the analytic continuation of $P_T$ over the path $\bs'(\tau)$.

\begin{theorem}[\cite{Purbhoo}]
 Let $T,T'\in\SYT(\maxp;\bs)$.
 Then $T_\bt$ is dual equivalent to $T'|_{\bt}$ if and only if 
 $P'_T(1)=P'_{T'}(1)$.
\end{theorem}

%
%
\subsection{Degree of the real Wronski map}
Recall that the complex Wronski map $\Wr\colon\Gr\to\PP^{(n+1)(d-n)}$ has
degree~\eqref{Eq:Grass_degree} 
\[
   \delta(\bi_{n,d})\ = \ [(n{+}1)(d{-}n)]!
          \frac{1!2!\dotsb n!}{(d{-}n)!(d{-}n{+}1)!\dotsb d!}\ .
\]
If we restrict the domain to the real Grassmannian, we get the 
real Wronski map $\Wr_\R\colon\Gr_\R\to\R\PP^{(n+1)(d-n)}$.
By Theorem~\ref{Th:strong}, over the locus of polynomials with $(n+1)(d-n)$ distinct real
roots, this is a $\delta(\bi_{n,d})$-to-one cover.
Eremenko and Gabrielov~\cite{EG02c} studied this real Wronski map, computing its topological
degree. 

This requires some explanation, for real projective spaces and Grassmannians are not always
orientable, and hence maps between them do not necessarily have a degree.
However, the Wronski map can be lifted to their orienting double covers, after which its degree
is well-defined up to a sign.
By the Pl\"ucker formula, the Wronski map restricted to the big Schubert cell $\Gr^\circ_\R$ of
the Grassmannian is a finite, proper map to $\R^{(n+1)(d-n)}$, realized as the space of monic
real polynomials of degree $(n+1)(d-n)$.
The compute the degree of the Wronski map over this big cell.

Fix a standard tableau $T_0\in\SYT(\maxp)$.
Given any tableau $T\in\SYT(\maxp)$, let $\sigma_T$ be the permutation in $\calS_{(n+1)(d-n)}$
with $\sigma_T(i)=j$ if the entries $i$ in $T_0$ and $j$ in $T$ occupy the same cell of
$\maxp$. 
Define
\[
   \DeCo{\delta(\bi_{n,d})_\R}\ :=\ 
     \sum_{T\in{\rm SYT}(\includegraphics[height=6pt]{pictures/max.eps})} |\sigma_T|\,,
\]
where $|\sigma|=\pm$ is the sign of the permutation $\sigma$.
  
\begin{theorem}[Theorem~2 of~\cite{EG02c}]
 $\deg \Wr_\R=\delta(\bi_{n,d})_\R$.
\end{theorem}

This statistic, $\delta(\bi_{n,d})_\R$, was computed by White~\cite{Wh}, who showed that it
vanishes unless $d$ is even, and in that case it equals
\[
   \frac{1!2!\dotsb(p{-}1)!(m{-}1)!(m{-}2)!\dotsb(m{-}p{+}1)!(\frac{mp}{2})!}
   {(m{-}p{+}2)!(m{-}p{+}4)!\dotsb(m{+}p{-}2)!\left(\frac{m-p+1}{2}\right)!%
     \left(\frac{m-p+3}{2}\right)!\dotsb\left(\frac{m+p-1}{2}\right)!}\ ,
\]
where $m:=\max\{n{+}1,d{-}n\}$ and $p:=\min\{n{+}1,d{-}n\}$.

The significance of these results is that $\delta(\bi_{n,d})_\R$ is a lower bound
for the number of real spaces of polynomials with given real Wronskian.
This gave the first example of a geometric problem possessing 
a nontrivial lower bound on its number of real solutions.
In the 1990's, Kontsevich~\cite{KM} determined the number \Blue{$N_d$} of
complex rational curves of degree $d$ interpolating $3d{-}1$ general points in the plane.
Work of Welschinger~\cite{W}, Mikhalkin~\cite{Mik}, and Itenberg, {\it et
  al.}~\cite{IKS03,IKS04} established a nontrivial lower bound \Blue{$W_d$} on the number of
real curves interpolating real points.
Not only is $W_d>0$, but
\[
   \lim_{d\to \infty} \frac{\log W_d}{\log N_d}\ =\ 1\qquad (!)
\]
More recently, Solomon~\cite{Sol} realized this number $W_d$ as the degree of a map.

Such lower bounds, if they were widespread, could have significant value for applications of
mathematics, as they are existence proofs for real solutions.
(On application of the nontriviality of $W_3=8$ is given in~\cite{FLT}.)
Initial steps in this direction were made in~\cite{SoSo,JW07}, which established lower bounds for
certain systems of sparse polynomials.

%
%
\section{Extensions of the Shapiro conjecture}\label{S:extensions}
The Shapiro conjecture for Grassmannians makes
sense for other flag manifolds. 
In this more general setting, it is known to fail, but in very interesting ways.
In some cases, we have been able to modify it to give a conjecture that holds under scrutiny.
The Shapiro conjecture also admits some appealing generalizations,
but its strongest and most subtle form remains open for Grassmannians.

%
%
\subsection{Lagrangian and Orthogonal Grassmannians}
The Lagrangian and orthogonal Grassmannians are closely related to the classical
Grassmannian.
For each of these, the Shapiro conjecture is particularly easy to state.

The (odd) orthogonal Grassmannian requires a nondegenerate symmetric
bilinear form $\langle \cdot,\cdot\rangle$ on $\C^{2n+1}$.
This vector space has a basis $e_1,\dotsc,e_{2n+1}$ such that 
\[
   \langle e_i,\, e_{2n+2-j} \rangle\ =\ \delta_{i,j}\,.
\]
A subspace $V$ of $\C^{2n+1}$ is \DeCo{{\sl isotropic}} if $\langle V,V\rangle=0$.
Isotropic subspaces have dimension at most $n$.
The \DeCo{{\sl (odd) orthogonal Grassmannian $OG(n)$}} is the set of all maximal 
($n$-dimensional) isotropic subspaces $V$ of $\C^{2n+1}$.
This variety has dimension $\binom{n{+}1}{2}$.

The Shapiro conjecture for $OG(n)$ begins with a particular rational normal curve $\gamma$
having parametrization
 \begin{multline*}
  \quad t\ \longmapsto\ e_1\ +\ t e_2\ +\ \frac{t^2}{2} e_3\ +\ \dotsb\ +\
      \frac{t^n}{n!}e_{n+1}\ \\ -\ \frac{t^{n+1}}{(n+1)!}e_{n+2}
     +\ \frac{t^{n+2}}{(n+2)!}e_{n+3}
     \ -\ \dotsb\ +\ (-1)^n\frac{t^{2n}}{(2n)!}e_{2n+1}\,. \quad
 \end{multline*}
This has special properties with respect to the form $\langle\cdot,\cdot\rangle$.
For $t\in\C$, define the flag $F_\bullet(t)$ in $\C^{2n+1}$ by
\[
   \DeCo{F_i(t)}\ :=\ \mbox{Span}\{ \gamma(t),\, \gamma'(t)\,,
          \dotsc,\, \gamma^{(i-1)}(t)\}\,.
\]
Then $F_\bullet(t)$ is \DeCo{{\sl isotropic}} in that
\[
    \langle F_i(t),\, F_{2n+1-i}(t)\rangle\ =\ 0\,.
\]
In general, an isotropic flag $F_\bullet$ of $\C^{2n+1}$ is a flag such that 
$\langle F_i, F_{2n+1-i}\rangle = 0$.

Schubert varieties of $OG(n)$ are defined with respect to an isotropic flag, $F_\bullet$,
and are the restriction of Schubert varieties of $\G(n{-}1,2n)$---the Grassmannian of $n$
dimensional subspaces of $\C^{2n+1}$---under the inclusion $OG(n)\hookrightarrow \G(n{-}1,2n)$.
Schubert varieties for $OG(n)$ are indexed by \DeCo{{\sl strict partitions}}, which are
integer sequences 
\[
   \bk\ \colon\ n\ \geq\ \kappa^1\ >\ \kappa^2\ >\ \dotsb\ >\ \kappa^k\ >\ 0\,.
\]
Set $\DeCo{\|\bk\|}=\kappa^1+\dotsb+\kappa^k$.
We do not give the precise relation between these indices and the ramification sequences
$\ba$ of Section~\ref{S:Shapiro}, but this may be done using the descriptions given
in~\cite[\S~6.1]{FuPr} or~\cite{So_iso}. 
Write \DeCo{$X_\sbk F_\bullet$} for the Schubert variety of $OG(n)$ defined by the Schubert
index $\bk$ and an isotropic flag $F_\bullet$.
It has codimension $\|\bk\|$.
A Schubert problem is a list $(\bk_1,\dotsc,\bk_m)$ of Schubert indices such that 
\[
   \|\bk_1\|+\|\bk_2\|+\dotsb+\|\bk_m\|\ =\ \dim OG(n)\ =\ \binom{n+1}{2}\,.
\]

We state the Shapiro conjecture for $OG(n)$.

\begin{conjecture}\label{Co:OG_n}
  If\/ $(\bk_1,\dotsc,\bk_m)$ is a Schubert problem for $OG(n)$ and 
  $s_1,\dotsc,s_m$ are distinct real numbers, then the intersection
\[
   X_{\sbk_1} F_\bullet(s_1)\ \bigcap\ 
   X_{\sbk_2} F_\bullet(s_2)\ \bigcap\ \dotsb\ \bigcap\ 
   X_{\sbk_m} F_\bullet(s_m)
\]
 is transverse with all points real.
\end{conjecture}

Besides optimism based upon the validity of the Shapiro conjecture for Grassmannians, the
evidence for Conjecture~\ref{Co:OG_n} comes in two forms.
Several tens of thousands of instances have been checked with a computer and  
when each $\|\bk_i\|=1$ and the points $s_i$ are clustered~\eqref{Eq:cluster}, the
intersection  is transverse with all points real~\cite{So00b}.\smallskip 

There is a similar story but with a different outcome for the Lagrangian Grassmannian.
Let $\langle \cdot,\cdot\rangle$ be a nondegenerate skew symmetric bilinear form on
$\C^{2n}$. 
This vector space has a basis $e_1,\dotsc,e_{2n}$ such that 
\[
   \langle e_i,\, e_{2n+1-j} \rangle\ =\ 
    \left\{\begin{array}{rcl} \delta_{i,j}&\ &\mbox{if $i\leq n$}\\
                             -\delta_{i,j}&\ &\mbox{if $i > n$}
      \end{array}\right. \ .
\]
Isotropic subspaces in $\C^{2n}$ may have any dimension up to $n$, and those of maximal
dimension are called \DeCo{{\sl Lagrangian subspaces}}.
The \DeCo{{\sl Lagrangian Grassmannian $LG(n)$}} is the set of all 
Lagrangian subspaces $V$ of $\C^{2n}$.
This variety has dimension $\binom{n{+}1}{2}$.

For the Shapiro conjecture for $LG(n)$, we have the rational normal curve $\gamma$
with parametrization
 \begin{multline*}
  \quad t\ \longmapsto\ e_1\ +\ t e_2\ +\ \frac{t^2}{2} e_3\ +\ \dotsb\ +\
  \frac{t^n}{n!}e_{n+1}\ \\-\ \frac{t^{n+1}}{(n+1)!}e_{n+2}
     +\ \frac{t^{n+2}}{(n+2)!}e_{n+3}
     \ -\ \dotsb\ +\ (-1)^{n-1}\frac{t^{2n-1}}{(2n-1)!}e_{2n}\,. \quad
 \end{multline*}
For $t\in\C$, define the flag $F_\bullet(t)$ in $\C^{2n+1}$ by
\[
   \DeCo{F_i(t)}\ :=\ \mbox{Span}\{ \gamma(t),\, \gamma'(t)\,,
          \dotsc,\, \gamma^{(i-1)}(t)\}\,.
\]
The flag $F_\bullet(t)$ is \DeCo{{\sl isotropic}} in that
\[
    \langle F_i(t),\, F_{2n-i}(t)\rangle\ =\ 0\,.
\]
More generally, an isotropic flag $F_\bullet$ of $\C^{2n}$ is a flag such that 
$\langle F_i, F_{2n-i}\rangle = 0$.

As with $OG(n)$, given an isotropic flag, Schubert varieties for $LG(n)$ are induced from
Schubert varieties of $\G(n{-}1,2n{-}1)$ by the inclusion
$LG(n)\hookrightarrow\G(n{-}1,2n{-}1)$. 
Schubert varieties \DeCo{$X_\sbk F_\bullet$} of $LG(n)$ are also indexed by strict partitions
$\bk$ and $\|\bk\|$ is the codimension of $X_\sbk F_\bullet$.
We give the relation between strict partitions for $LG(n)$ and ramification sequences
for $\G(n{-}1,2n{-}1)$.
Given a strict partition $\bk\colon n\geq \kappa^1>\dotsb>\kappa^k$, let
$\bmu\colon 0<\mu_1<\dotsb<\mu_{n-k}$ be the complement of the set 
$\{\kappa^1,\dotsc,\kappa^k\}$ in $\{1,2,\dotsc,n\}$. 
Call $k$ the \DeCo{{\sl length}} of the strict partition $\bk$.
For example, if $n=6$ and $\bk=4,2$, then $k=2$ and $\bmu=1,3,5,6$.
If we define $\DeCo{\ba(\bk)}=(a_0,\dotsc,a_{n-1})$ to be the sequence
\[
   0\leq n-\kappa^1<\dotsb<n-\kappa^k\;<\; n{-}1+\mu_1 <\dotsb< n{-}1+\mu_{n-k}\leq 2n{-}1\,, 
\]
then $X_\sbk F_\bullet=\Omega_{\ba(\sbk)}F_\bullet\cap LG(n)$, so that
\[
   X_\sbk F_\bullet\ =\ \{ V\in LG(n)\mid F_{2n-a_j}\geq n-j,\ 
          \mbox{for}\ j=0,1,\dotsc,n{-}1\}\,.
\]

A Schubert problem is a list $(\bk_1,\dotsc,\bk_m)$ such that 
\[
   \|\bk_1\|+\|\bk_2\|+\dotsb+\|\bk_m\|\ =\ \dim LG(n)\ =\ \binom{n+1}{2}\,.
\]
The obvious generalization of Theorem~\ref{Th:MTV_1} and Conjecture~\ref{Co:OG_n} to
$LG(n)$ turns out to be false.
We offer a modification that we believe is true.
Belkale and Kumar~\cite{BK} define a notion they call Levi movability.
A Schubert problem $(\bk_1,\dotsc,\bk_m)$ for $LG(n)$ is \DeCo{{\sl Levi movable}} if
the corresponding Schubert indices, $(\ba(\bk_1),\dotsc,\ba(\bk_m))$ also form a Schubert
problem for $\G(n{-}1,2n{-}1)$.
Unraveling the definitions shows that this is equivalent to having the lengths of
the strict partitions  $(\bk_1,\dotsc,\bk_m)$ sum to $n$.

\begin{conjecture}\label{Co:LG_n}
  If $(\bk_1,\dotsc,\bk_m)$ is a Schubert problem for $LG(n)$ and 
  $s_1,\dotsc,s_m$ are distinct real numbers, then the intersection
\[
   X_{\sbk_1} F_\bullet(s_1)\  \bigcap\ 
   X_{\sbk_2} F_\bullet(s_2)\  \bigcap\ \dotsb\  \bigcap\ 
   X_{\sbk_m} F_\bullet(s_m)
\]
 is transverse.
 If $(\bk_1,\dotsc,\bk_m)$ is Levi movable, then all points of intersection are real,
 but if it is not Levi movable, then no point in the intersection is real.
\end{conjecture}

The strongest evidence in favor of Conjecture~\ref{Co:LG_n} is that it is true when 
the Schubert problem $(\bk_1,\dotsc,\bk_m)$ is Levi movable.
This follows from the definition of Levi movable and the Shapiro conjecture for Grassmannians.
Further evidence is that if each $\bk_i$ is simple in that $\|\bk_i\|=1$, then 
a local version, similar to Theorem~\ref{Th:local} but without transversality, is true~\cite{So00b}.
That is, if the $s_i$ are clustered~\eqref{Eq:cluster}, then no point in the intersection is 
real.
Lastly, several tens of thousands of instances have been checked with a computer.

%
%
\subsection{Monotone conjecture for flag manifolds}
The Shapiro conjecture was originally made for the
classical (type-$A$) flag manifold, where is fails spectacularly.
It is false for the first nontrivial Schubert problem on
a flag variety that is not a Grassmannian.
Namely, the geometric problem of partial flags consisting of a line $\ell$ lying on a plane 
$\Lambda$ in 3-dimensional
space where $\ell$  meets three fixed lines and $\Lambda$ 
contains two fixed points. 
 
This is just the problem of four lines in disguise.
Suppose that $p$ and $q$ are the two fixed points that $\Lambda$ is required to contain.
Then $\Lambda$ contains the line $\overline{p,q}$ they span.
Since $\ell\subset\Lambda$, it must meet $\overline{p,q}$.
As $\ell$ must also meet three lines, this problem reduces to
the problem of four lines.
In this way, there are two solutions to this Schubert problem.

Now let us investigate the Shapiro conjecture for this Schubert problem, which
posits that both flags $\ell\subset\Lambda$ will be real, if we require that $\ell$
meets three fixed tangent lines to a rational curve and $\Lambda$ contains two fixed
points of the rational curve.
Let $\gamma$ be the rational normal curve~(1) from the Introduction and suppose that 
the three fixed lines of our problem are its tangent lines
$\ell(-1)$, $\ell(0)$, and $\ell(1)$.
These lines lie on the hyperboloid $H$ of one sheet~(2).
Here is another view of these lines, the curve $\gamma$, and the hyperboloid.
 \[
  \begin{picture}(260,163)(0,-1)
   \put(0,0){\includegraphics[height=160pt]{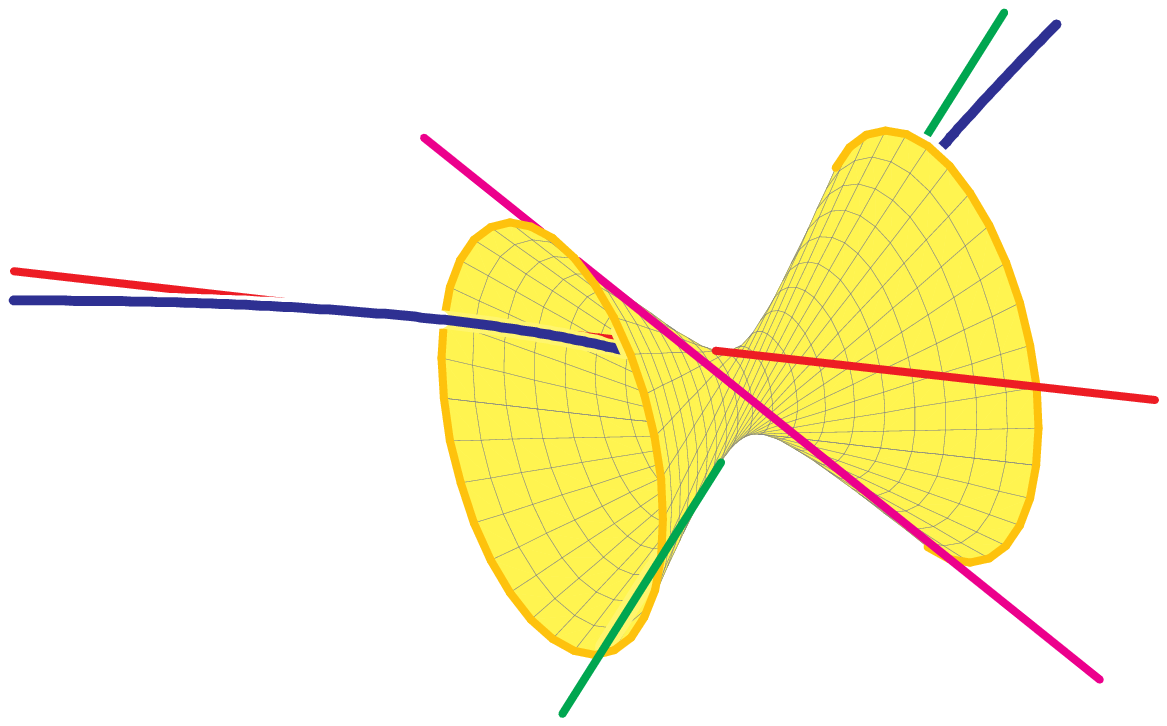}}
   \put(100,3){$\ell(-1)$} \put(242,17){$\ell(0)$}
   \put(1,105){$\ell(1)$} \put(1,83){$\gamma$}
   \put(225,115.2){$H$}
  \end{picture}
 \]
If we require $\ell$ to meet the three tangent lines $\ell(-1)$, $\ell(0)$, and
$\ell(1)$ and $\Lambda$ to contain the two points $\gamma(v)$ and $\gamma(w)$ of
$\gamma$, then $\ell$ also meets the line $\lambda(v,w)$ spanned by these two points.
As in the Introduction, the lines $\ell$ that we seek will come from points where the
secant line $\lambda(v,w)$ meets $H$.

\begin{figure}[htb]
 \[
  \begin{picture}(240,130)(-18,-10)
   \put(0,0){\includegraphics[height=110pt]{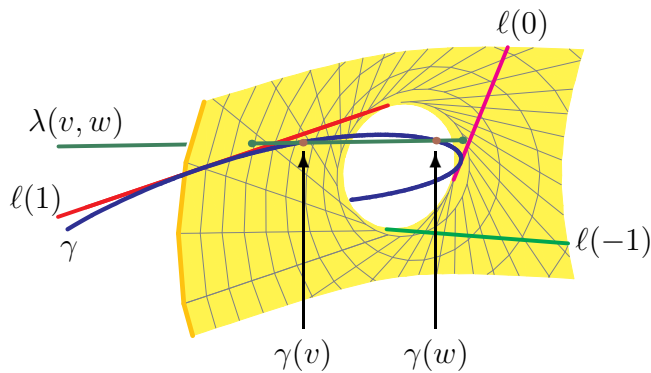} }
   \put(-18,49){$\ell(1)$} \put(196,33){$\ell(-1)$} 
   \put(166,115){$\ell(0)$}
   \put(1,31){$\gamma$} \put(-11,78){$\lambda(v,w)$}
  \thicklines
    \put( 82,-10){$\gamma(v)$}  \put( 93,3){\vector(0,1){67}}
    \put(131,-10){$\gamma(w)$}  \put(143,3){\vector(0,1){67}}
  \end{picture}
\]
\caption{A secant line meeting $H$.}\label{F:throat}
\end{figure}

Figure~\ref{F:throat} shows an expanded view down the throat of the hyperboloid, 
with a secant line $\lambda(v,w)$ that meets the hyperboloid in two points. 
For these points $\gamma(v)$ and $\gamma(w)$ there will be two real flags
$\ell\subset\Lambda$ satisfying our conditions.
This is consistent with the Shapiro conjecture.

\begin{figure}[htb]
 \[
  \begin{picture}(335,220)
   \put(0,0){\includegraphics[height=220pt]{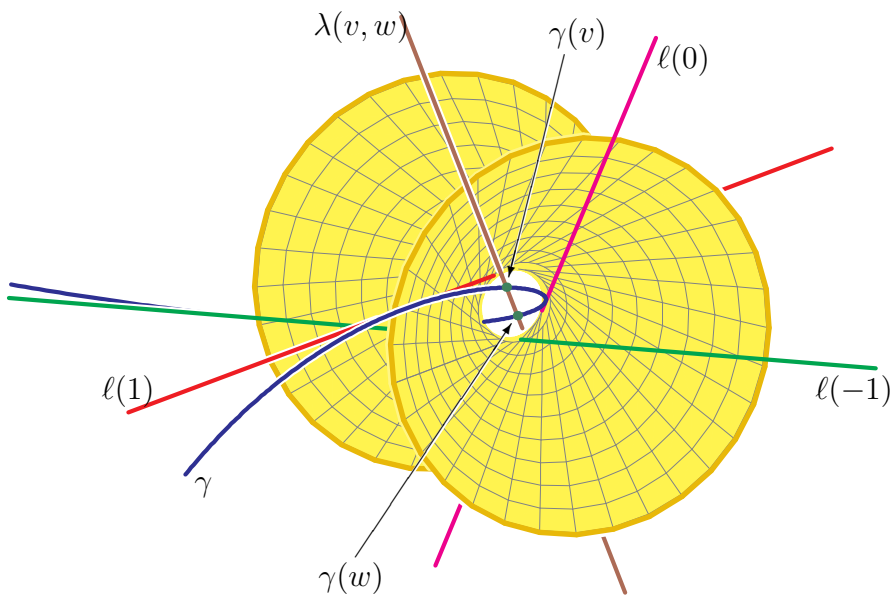}}
   \put(35,75){$\ell(1)$} \put(305,75){$\ell(-1)$} 
   \put(245,200){$\ell(0)$}
   \put(70,40){$\gamma$}\put(115,212.5){$\lambda(v,w)$}
   \put(204,210){$\gamma(v)$}  
       \put(210,205){\vector(-1,-4){21}}
   \put(117,4){$\gamma(w)$}  
       \put(131,15){\vector(2,3){59}} 
  \end{picture}
 \]
\caption{A secant line not meeting $H$.}\label{F:two}
\end{figure}

In contrast, Figure~\ref{F:two} shows a secant line $\lambda(v,w)$ that does not meet the
hyperboloid in any real points.
For these points $\gamma(v)$ and $\gamma(w)$, neither flag $\ell\subset\Lambda$ satisfying 
our conditions is real.
This is a counterexample to the Shapiro conjecture.\smallskip

This failure of the Shapiro conjecture is however quite interesting.
If we label the points $-1,0,1$ with 1 (conditions on the line) and $v,w$ by 2
(conditions on the plane), then along $\gamma$ they occur in order 
 \begin{equation}\label{Eq:condition_order}
   \mbox{$11122$ in Figure~\ref{F:throat} \ \ and \ \ $11212$ in Figure~\ref{F:two}.}
 \end{equation}
The sequence for Figure~\ref{F:throat} is \DeCo{{\it monotone}} increasing and in this case 
both solutions are always real, but the sequence for Figure~\ref{F:two} is not monotone.
This example suggests a way to correct the Shapiro conjecture, that we call 
the \DeCo{{\sl monotone conjecture}}.

Specifically, let $\bn\colon 0\leq n_1<\dotsb<n_k<d$ be a sequence of integers.
The manifold $\Fl$ of flags of type $\bn$ is the set of all sequences of subspaces
\[
  E_\bullet\ \colon\ E_{n_1}\ \subset\ E_{n_2}\ \subset\ 
    \dotsb\ \subset\ E_{n_k}\ \subset\ \C_d[t]
\]
with $\dim E_{n_i}=n_i+1$.
The forgetful map $E_\bullet\mapsto E_{n_i}$ induces a projection
\[
  \DeCo{\pi_i}\ \colon\ \Fl\ \longrightarrow\ \G(n_i,d)
\]
to a Grassmannian.
A \DeCo{{\sl Grassmannian Schubert variety}} is a subvariety of $\Fl$ of the form
$\pi_i^{-1} \Omega_{\ba} F_\bullet$.
That is, it is the inverse image of a Schubert variety in a Grassmannian projection.
Write \DeCo{$X_{(\ba,n_i)}F_\bullet$} for this Grassmannian Schubert variety and call
\DeCo{$(\ba,n_i)$} a Grassmannian Schubert condition.

A \DeCo{{\sl Grassmannian Schubert problem}} is a list
 \begin{equation}\label{Eq:GSP}
   (\ba^{(1)},n^{(1)}),\ 
   (\ba^{(2)},n^{(2)}),\ \dotsc,\ 
   (\ba^{(m)},n^{(m)}),\ 
 \end{equation}
of Grassmannian Schubert conditions satisfying
$|\ba^{(1)}|+\dotsb+|\ba^{(m)}|=\dim \Fl$.
We state the monotone conjecture.

\begin{conjecture}
 Let $\bigl((\ba^{(1)},n^{(1)}),\dotsc,(\ba^{(m)},n^{(m)})\bigr)$ be a Grassmannian
 Schubert problem for the flag variety $\Fl$ with
 $n^{(1)} \leq n^{(2)} \leq \dotsb \leq n^{(m)}$.
 Whenever $s_1<s_2<\dotsb<s_m$ are real numbers, the intersection
\[
   X_{(\ba^{(1)},n^{(1)})}F_\bullet(s_1)\ \bigcap\ 
   X_{(\ba^{(2)},n^{(2)})}F_\bullet(s_2)\ \bigcap\ \dotsb \ \bigcap\ 
   X_{(\ba^{(m)},n^{(m)})}F_\bullet(s_m)\,,
\]
 is transverse with all points of intersection real (when it is nonempty).
\end{conjecture}

There is significant evidence for this monotone conjecture.
First, the Shapiro conjecture for Grassmannians is the special case case 
when $m=1$ so then $\bn=n_1$ and $\Fl=\G(n_1,d)$:
the monotonicity condition $s_1<\dotsb<s_m$ is empty as any reordering of the 
Schubert conditions remains sorted. 

This conjecture was formulated in~\cite{RSSS}.
That project was based upon computer experimentation using  15.76 gigaHertz-years of
computing to study over 520 million instances of 1126 different Schubert problems on 29 flag
manifolds.
Some of this computation studied intersections of Schubert varieties that were not
necessarily monotone.
For example, consider the Schubert problem on $\Flv_{1<2,5}$,
 \begin{equation}\label{Eq:12-flag}
   \DeCo{( 0{<}2\,,\;1)^4}\,,\  \Mulberry{(0{<}1{<}3\,,\;2)^4}\,,
 \end{equation}
where the exponent indicates a repeated condition.
Table~\ref{table:12-flag} displays the computation on this Schubert problem.
 \begin{table}[htb]
  \begin{tabular}
   {|c||r|r|r|r|r|r|r|}\hline
   {} & \multicolumn{7}{c|}{Number of Real Solutions\rule{0pt}{11pt}}\\
   \cline{2-8}
        &0&2&4&6&8&10&12\rule{0pt}{13pt}\\\hline\hline
   \DeCo{1111}\Mulberry{2222}
     &  &   &   &   &   &   & 400000\\\hline
   \DeCo{11}\Mulberry{2}\DeCo{11}\Mulberry{222}   
     &  &   & 118 & 65425 & 132241 & 117504 & 84712\\ \hline
   \DeCo{111}\Mulberry{22}\DeCo{1}\Mulberry{22}
     &  &   & 104 & 65461 & 134417 & 117535 & 82483 \\ \hline
   \DeCo{11}\Mulberry{22}\DeCo{11}\Mulberry{22}
     &  &   & 1618 & 57236 & 188393 & 92580 & 60173  \\ \hline
   \DeCo{11}\Mulberry{2}\DeCo{1}\Mulberry{22}\DeCo{1}\Mulberry{2}
     &  &   & 25398 & 90784 & 143394 & 107108 & 33316 \\ \hline
   \DeCo{11}\Mulberry{22}\DeCo{1}\Mulberry{2}\DeCo{1}\Mulberry{2}
     &  & 2085 & 79317 & 111448 & 121589 & 60333 & 25228\\ \hline
   \DeCo{111}\Mulberry{2}\DeCo{1}\Mulberry{222}
     &  & 7818 & 34389 & 58098 & 101334 & 81724 & 116637 \\ \hline
   \DeCo{1}\Mulberry{2}\DeCo{1}\Mulberry{2}\DeCo{1}\Mulberry{2}\DeCo{1}\Mulberry{2}
     & 15923 & 41929 & 131054 & 86894 & 81823 & 30578 & 11799 \\ \hline
  \end{tabular}\vspace{5pt}

  \caption{The Schubert problem~\eqref{Eq:12-flag} on $\Flv_{\DeCo{1}<\Mulberry{2},5}$.}  
  \label{table:12-flag}
 \end{table}
%
%
The rows are labeled by different orderings of the conditions along the
rational normal curve $\gamma$ in the notation of~\eqref{Eq:condition_order}.
Each cell contains the number of computed instances with a given ordering
and number of real solutions.
The empty cells indicate no observed instances. 
Only the first row tests the monotone conjecture:
Each of the 400,000 computed instances had all 12 solutions real.
The other rows reveal a very interesting pattern; 
for nonmonotone orderings of the conditions along $\gamma$, not all solutions are always real
and there seems to be a lower bound on the number of real solutions.
Only in the last row, which represents the maximal possible
intertwining of the conditions, were  no real solutions observed.

A third piece of evidence for the monotone conjecture was provided by Eremenko,
et.~al~\cite{EGSV}, who showed that it is true for two-step flag manifolds, when
$\bn=d{-}2<d{-}1$. 
This is a special case of their main theorem, which asserts the reality of a 
rational function $\varphi$ with prescribed critical points on $\R\PP^1$ and certain prescribed  
coincidences $\varphi(v)=\varphi(w)$, when $v,w$ are real.

This result of Eremenko, et.~al can be described in terms of $\G(d{-}2,d)$, where it
becomes a statement about real points in an intersection of Schubert varieties 
given by flags that are secant to a rational normal curve in a particular way.
This condition on secant flags makes sense for any Grassmannian, and the resulting
secant conjecture is also a generalization of the Shapiro conjecture.

A flag $F_\bullet$ is \DeCo{{\sl secant along an arc $I$}} of a rational normal curve $\gamma$ 
if every subspace in the flag is spanned by its intersections with $I$.
A collection of flags that are secant to $\gamma$ is \DeCo{{\sl disjoint}} if they are
secant along disjoint arcs of $\gamma$.
The \DeCo{{\sl secant conjecture}} asserts that a Schubert problem given by 
disjoint secant flags has all solutions real.
We give a more precise statement.

\begin{conjecture}
  If $(\ba_1,\dotsc,\ba_m)$ is a Schubert problem for $\Gr$ and 
  $F_\bullet^1,\dotsc,F_\bullet^m$ are disjoint secant flags, then the intersection
\[
    \Omega_{\ba_1} F_\bullet^1\ \bigcap\ 
    \Omega_{\ba_2} F_\bullet^2\ \bigcap\ \dotsb \ \bigcap\ 
    \Omega_{\ba_m} F_\bullet^m
\]
 is transverse with all points real.
\end{conjecture}

The main result of~\cite{EG02a} is that an intersection of Schubert varieties in 
$\G(d-2,d)$ given by disjoint secant flags is transverse with all points real.
The Shapiro conjecture is a limiting case of the secant conjecture, as the flag osculating
$\gamma$ at a point $s$ is the limit of flags that are secant along arcs that shrink to the
point $s$.

\begin{figure}[htb]
\[
  \begin{picture}(370,153)(0,-5)
   \put(0,7){\includegraphics[height=140pt,viewport=5 78 430 240,clip]{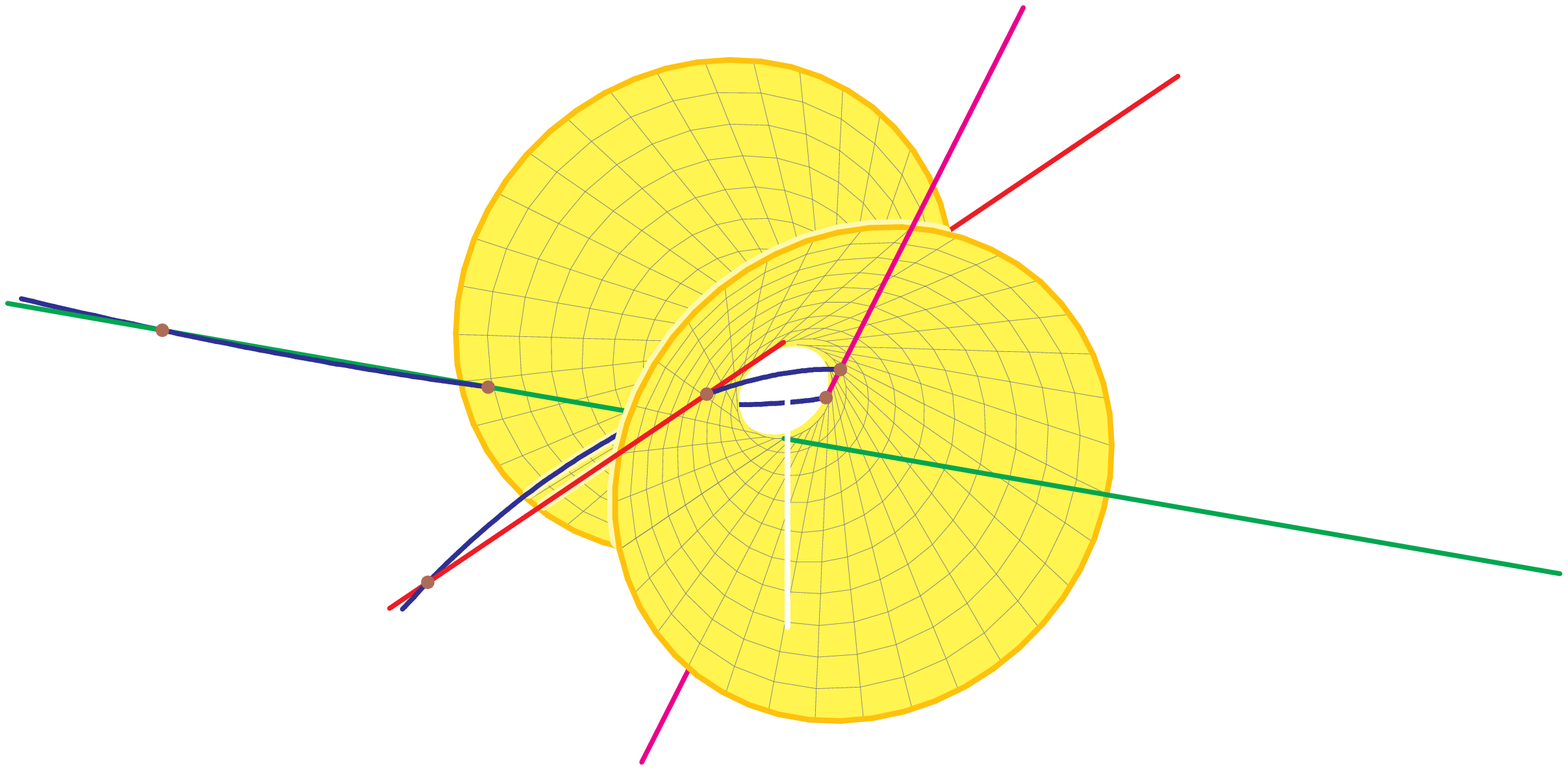}}
   \thicklines
   \put(139,49){$\gamma$}\put(148,51){\vector(2,-1){25}}
   \put(294.14,-2){\vector(0,1){97}}
   \put(290,-12){$I$}
  \end{picture}
\]
\caption{The problem of four secant lines.}
\label{F:secant}
\end{figure}

Consider this secant conjecture for the problem of four lines.
The hyperboloid in Figure~\ref{F:secant} contains three lines that are secant to 
$\gamma$ along disjoint arcs.
Any line secant along the arc {$I$} (which is disjoint from the other three arcs)
meets the hyperboloid in two points, giving two real solutions to this
instance of the secant conjecture.


This secant conjecture is currently being studied on a supercomputer whose day job is
calculus instruction.
For each of hundreds of Schubert problems, thousands to millions of instances of the
secant conjecture are being tested, and much more.
The \DeCo{{\sl overlap number}} measures how far a collection of secant flags is from
being disjoint, and it is zero if and only if the flags are disjoint.
This experiment tests instances of the secant conjecture {\sl and} near
misses when the flags have low overlap number.
The results (number of real solutions vs.~overlap number) are stored in a publicly
accessible database accessible from the webpage~\cite{secant}.
In the first nine months of operation, this has studied over 1.3 billion instances of
Schubert problems and consumed over 600 gigaHertz-years of computing.

Table~\ref{T:Sample_table} shows the results
for a Schubert problem  with 16 solutions on $\G(2,5)$.
Computing the $20,000,000$ instances of this problem used 4.473 gigaHertz-years.
The rows are labeled with the even integers from 0 to 16 as the number of real solutions has
the same parity as the number of complex solutions. 
\begin{table}[htb]
%

\noindent{\small 

\noindent\begin{tabular}{|r||r|r|r|r|r|r|r|c||r|}\hline
  \multirow{2}{29pt}{\# real solns.} & \multicolumn{9}{c|}{Overlap Number\qquad\rule{0pt}{11pt}}\\
   \cline{2-10}
   & 0&1&2&3&4&5&6&\ $\dotsb$\ &Total\\\hline\hline
0 &&&&&&&20&\ $\dotsb$\ & 7977\\\hline
2 &&&&&&&116&\ $\dotsb$\ & 88578\\\hline
4 &&&&6154&23561&526&3011&\ $\dotsb$\ & 542521\\\hline
6 &&&&25526&63265&2040&9460&\ $\dotsb$\ & 1571582\\\hline
8 &&&&33736&78559&2995&13650&\ $\dotsb$\ & 2834459\\\hline
10 &&&&25953&39252&2540&11179&\ $\dotsb$\ & 3351159\\\hline
12 &&&&35578&44840&3271&14160&\ $\dotsb$\ & 2944091\\\hline
14 &&&&17367&17180&1705&7821&\ $\dotsb$\ & 1602251\\\hline
16 &4568553&&182668&583007&468506&36983&83169&\ $\dotsb$\ &7057382 \\\hline\hline
Total &4568553&&182668&727321&735163&50060&142586&\ $\dotsb$\ &20000000 \\\hline
\end{tabular}

}\vspace{4pt}
\caption{Number of Real solutions v.s. overlap number.}
\label{T:Sample_table}
\end{table}
The column with overlap number 0 represents tests of the secant conjecture.
Since its only entries are in the row for 16 real solutions, the secant
conjecture was verified in $4,568,553$ instances.
The column labeled 1 is empty because flags for this
problem cannot have overlap number 1.
The most interesting feature is that for overlap number 2, all
solutions were still real, while for overlap numbers 3, 4, and 5, at least 4 of the 16
solutions were real, and only with overlap number 6 and greater does the Schubert problem have
no real solutions. 
This inner border, which indicates that the secant conjecture does not completely fail when
there is small overlap, is found on many of the other problems that we investigated and is a new
phenomenon that we do not understand.
A description of the technical aspects of this running experiment is given
in~\cite{secant_experiment}.

\subsection{Discriminant conjecture}
Despite the proofs of Theorems~\ref{Th:MTV_1}
and~\ref{Th:strong} (weak and strong form of the Shapiro for Grassmannians), the strongest and
most subtle form of that conjecture remains open.

The discriminant of a polynomial $W=\prod_i(t-s_i)$ is
$\prod_{i<j}(s_i-s_j)^2$, the symmetric function of its roots having lowest degree
that vanishes when $W$ has a double root.
More generally, suppose that we have a family of polynomial systems in a space $X$ that are
parametrized by a space $S$.
(For example, the intersection in Theorem~\ref{Th:strong},
\[
    \Omega_{\ba^{(1)}} F_\bullet(s_1)\ \bigcap\ 
    \Omega_{\ba^{(2)}} F_\bullet(s_2)\ \bigcap\ \dotsb\ \bigcap\ 
    \Omega_{\ba^{(m)}} F_\bullet(s_m) \leqno{(\ref{Eq:tr_int})}
\]
in which $X=\Gr$ and $S$ is $\C^m$ or $(\PP^1)^m$.)
Then the discriminant variety of this system is the subvariety $\Sigma\subset S$ 
where the system is not transverse.
This is expected to be a hypersurface, and the \DeCo{{\sl discriminant}} of the
system is the function that defines $\Sigma$.

By Theorem~\ref{Th:strong} this discriminant does not vanish when the parameters $s_i$ are real
and distinct.
However, in the few cases when it has been computed much more is true,
it is a sum of squares~\cite{So00} and therefore nonnegative.
For example, for the Schubert problem $\bi_{1,4}$ with 5 solutions, if we fix $s_5=0$ and
$s_6=\infty$, then the discriminant is a homogeneous polynomial of degree 20 in the four
variables $s_1,\dotsc,s_4$ with 711 terms, which turns out to be a sum of squares.
This is remarkable because Hilbert~\cite{Hi1888} showed that, except for $m=3$ and
$\deg=4$, not all nonnegative homogeneous polynomials in $m>2$ variables of degree more than 2
are sums of squares. 
Work of Blekherman~\cite{Bl} suggests that it is extremely rare for a nonnegative
polynomial to be a sum of squares.

\begin{conjecture}[Question 4 of~\cite{So00}]
 The discriminant of an intersection~$\eqref{Eq:tr_int}$ of a Schubert problem on a Grassmannian
 given by osculating flags is sum of squares in the parameters $s_1,\dotsc,s_m$.
\end{conjecture}

We conjecture that this remains true for any cominuscule flag variety, which includes the 
Lagrangian Grassmannian, the orthogonal Grassmannian, quadrics, as well as the two exceptional
cases $E_6/D_5$ and $E_7/E_6$.
There is also a form of this conjecture, Conjecture 2.10 of~\cite{RSSS}, involving preorders
for the semialgebraic set of monotone parameters $s_1<s_2<\dotsb<s_m$.

We close with the remark that we have not yet investigated the Shapiro conjecture for other
flag manifolds, and do not yet know when it fails, or how to repair the failures.
Also, the methods of Mukhin, Tarasov, and Varchenko only work for the Grassmannian, and it is
completely unclear how to even approach a proof of these generalizations.

\providecommand{\bysame}{\leavevmode\hbox to3em{\hrulefill}\thinspace}
\providecommand{\MR}{\relax\ifhmode\unskip\space\fi MR }
\providecommand{\MRhref}[2]{%
  \href{http://www.ams.org/mathscinet-getitem?mr=#1}{#2}
}
\providecommand{\href}[2]{#2}


\end{document}